\documentclass[reqno]{amsproc}
\usepackage[latin1]{inputenc}
\usepackage{amssymb}
\usepackage{hyperref}
\usepackage{amsmath}
\usepackage{latexsym}
\usepackage{cite}
\usepackage{amssymb}
\usepackage{amsfonts}
\usepackage{amscd}
\usepackage{xcolor}
\usepackage{amstext,amsmath,amssymb,amsfonts}
\newtheorem{theorem}{Theorem}

\newtheorem{corollary}[theorem]{Corollary}

\newtheorem{definition}[theorem]{Definition}

\newtheorem{proposition}[theorem]{Proposition}

\newcommand{\K}{\mathbb {K}}

\newcommand{\A}{\mathcal{A}}

\newcommand{\beq}{\begin{eqnarray}}
\newcommand{\eeq}{\end{eqnarray}}
\newcommand{\beqs}{\begin{eqnarray*}}
\newcommand{\eeqs}{\end{eqnarray*}}
\newcommand{\bpro}{\begin{pro}}
\newcommand{\epro}{\end{pro}}
\newcommand{\blem}{\begin{lem}}
\newcommand{\elem}{\end{lem}}
\newcommand{\bdfn}{\begin{dfn}}
\newcommand{\edfn}{\end{dfn}}
\newcommand{\bcor}{\begin{cor}}
\newcommand{\ecor}{\end{cor}}
\newcommand{\bthm}{\begin{thm}}
\newcommand{\ethm}{\end{thm}}
\newcommand{\bex}{\begin{ex}}
\newcommand{\eex}{\end{ex}}
\newcommand{\brmk}{\begin{rmk}}
\newcommand{\ermk}{\end{rmk}}
\newcommand{\bpr}{\begin{pr}}
\newcommand{\epr}{\end{pr}}
\newcommand{\benum}{\begin{enumerate}}
\newcommand{\eenum}{\end{enumerate}}
\newcommand{\bitem}{\begin{itemize}}
\newcommand{\eitem}{\end{itemize}}

\chardef\bslash=`\\
\numberwithin{equation}{section}
\numberwithin{table}{section}
\numberwithin{theorem}{section}
\DeclareMathOperator{\id}{id}

\title[ Yang-Baxter equation, $D-$equation, Frobenius algebras and Connes cocycles ]{Solutions of associative Yang-Baxter equation and $D-$equation in low dimensions and associated Frobenius algebras and Connes cocycles\footnote{Preprint: ICMPA-MPA/2015/08 } }

\author{Mahouton Norbert Hounkonnou$^\ast$}
\address[$\ast$]{University of Abomey-Calavi,
International Chair in Mathematical Physics and Applications,
ICMPA-UNESCO Chair, 072 BP 50, Cotonou, Rep. of Benin}
\email{norbert.hounkonnou@cipma.uac.bj, with copy to hounkonnou@yahoo.fr}
\author{Gb\^ev\`ewou Damien  Houndedji$^\dagger$}
\address[$\dagger$]{University of Abomey-Calavi,
International Chair in Mathematical Physics and Applications,
ICMPA-UNESCO Chair, 072 BP 50, Cotonou, Rep. of Benin}
\email{ houndedjid@gmail.com}
\begin{document}
\maketitle

\today

\bigskip
\begin{abstract}
This work addresses some relevant characteristics of associative algebras in low dimensions.
Especially, given $ 1 $ and $ 2 $ dimensional associative algebras, we explicitly solve
associative Yang-Baxter equations and use skew-symmetric solutions to perform double constructions
of Frobenius algebras. Besides, we determine related compatible dendriform algebras and solutions of their $ D- $equations. Finally, using symmetric solutions of the latter equations, we proceed to double constructions of corresponding Connes cocycles.  
\end{abstract}
\section{Introduction}

 A (symmetric) Frobenius algebra which is an associative
algebra with a (symmetric) non-degenerate invariant bilinear form
is an important object in both mathematics and mathematical
physics. It plays a key role in the study of many topics, such as
statistical models over 2-dimensional graphs \cite{[MB]} and topological quantum field theory
\cite{[JK]}. On the other hand, a
non-degenerate  Connes cocycle is an associative algebra
with a non-degenerate antisymmetric bilinear form being a cyclic
1-cocyle in the sense of Connes \cite{[AC]}. It corresponds to the original definition of
cyclic cohomology by Connes and hence is important in the study of
noncommutative geometry.

However, it is not easy to construct Frobenius algebras or
non-degenerate Connes cocycles explicitly, that is, both the
explicit examples of  Frobenius algebras and non-degenerate Connes
cocycles are lacked. In \cite{[C.Bai5]}, some special
constructions (namely, double constructions) of both two objects
were given in terms of bialgebra structures and certain algebraic
equations. In particular, it provides an approach to construct
both (symmetric) Frobenius algebras and non-degenerate Connes
cocycles from solving certain algebraic equations. Explicitly, a
(symmetric) Frobenius algebra can be obtained from an
anti-symmetric solution of associative Yang-Baxter equation in an
associative algebra, whereas a non-degenerate Connes cocycle can
be obtained from a symmetric solution of $D$-equation in a
dendriform algebra which is the underlying algebraic structure of
a non-degenerate Connes cocycle. Note that both associative
Yang-Baxter equation and dendriform algebras appeared more early
in some other fields. For example, the associative Yang-Baxter
equation was introduced by Aguiar \cite{[M1]} to study  the cases
of principal derivations for the infinitesimal bialgebras given by
Joni and Rota \cite{[Joni]} to provide an algebraic framework for
the calculus of divided difference, whereas dendriform algebras
were introduced by Loday \cite{[N8]} with motivation from
algebraic K-theory and were studied quite extensively with
connections to several areas in mathematics and physics, like
operads \cite{[N10]}, homology \cite{[N4]}, \cite{[N5]},
arithmetics \cite{[N9]}  and quantum field theory \cite{[N3]}.

In this paper, under the above framework, we will give the
explicit study in low dimensions. The paper is organized as
follows. In Section 2, we give some basic notions and results on
the double constructions of Frobenius algebras and Connes
cocycles. In Section 3, we give the explicit study in dimension 1.
Specifically, given $ 1-$dimensional associative algebras, we
explicitly solve
  associative Yang-Baxter equations and use skew-symmetric solutions to perform double constructions
   of Frobenius algebras. Then, we determine related compatible dendriform algebras and solutions of their $D$-equations.
  Finally, using symmetric solutions of the latter equations, we proceed to double constructions of corresponding Connes cocycles. In Section 4, similar calculations are performed in dimension 2.
In Section 5, we give some  concluding remarks.

 \section{Preliminaries}
In this section, we give a quick overview on main definitions and  fundamental results essentially known from \cite{[C.Bai5]}. See also  \cite{[M1]},  \cite{[N14]}-\cite{[N15]}, \cite{[N8]} and  \cite{[N6]} and the references therein.
\begin{definition} \label{d1}
A bilinear form $ \mathcal{B}(\cdot, \cdot ) $ on an associative algebra $ \mathcal{A} $ is \textbf{invariant} if
\beqs
\mathcal{B}(xy,z) = \mathcal{B}(x, yz) \mbox { for all } x, y, z \in  \mathcal{A} .
\eeqs
\end{definition}

\begin{definition} \label{d2}
An antisymmetric bilinear form $ \omega(\cdot, \cdot)  $ on an associative algebra $ \mathcal{A} $ is a \textbf{cyclic 1-cocycle in the sense of Connes} if
\beq \label{eq4}
\omega(xy,z) + \omega(yz,x) + \omega(zx,y) = 0 \mbox { for all } x, y, z \in \mathcal{A}.
\eeq
For simplicity,
 $ \omega $ is called a \textbf{Connes cocycle}.
\end{definition}

\begin{definition} \label{d3}
A Frobenius algebra $ (\mathcal{A}, \mathcal{B}) $ is an associative algebra $ \mathcal{A} $ with a non-degenerate invariant bilinear form $ \mathcal{B} (\cdot, \cdot) $.
It is symmetric if $  \mathcal{B} $ is symmetric.
\end{definition}
\begin{definition}  \label{d5}
We call $(\mathcal{A}, \mathcal{B})$ a   double
construction of a (symmetric) Frobenius algebra associated to
$\mathcal{A}_1$ and ${\mathcal A}_1^*$ if it satisfies the
conditions \benum \item[(1)] $ \mathcal{A} = \mathcal{A}_{1}
\oplus \mathcal{A}^{\ast}_{1} $ as the direct sum of vector
spaces; \item[(2)] $ \mathcal{A}_{1} $ and $
\mathcal{A}^{\ast}_{1} $ are associative subalgebras of $
\mathcal{A} $; \item[(3)] $ \mathcal{B} $ is the natural symmetric
bilinear form on $ \mathcal{A}_{1} \oplus \mathcal{A}^{\ast}_{1} $
given by \beq \label{eq2}
 \mathcal{B}(x + a^{\ast}, y + b^{\ast}) = \langle x, b^{\ast} \rangle +  \langle a^{\ast}, y \rangle \mbox { for all } x, y \in \mathcal{A}_{1}, a^{\ast}, b^{\ast} \in \mathcal{A}^{\ast}_{1},
\eeq
where $ \langle  , \rangle $ is the natural pair between the vector space $ \mathcal{A}_{1} $ and its dual space  $ \mathcal{A}^{\ast}_{1} $.
\eenum

We call $ (\mathcal{A}, \omega) $ 
a Connes cocycle associated to
$\mathcal{A}_1$ and ${\mathcal A}_1^*$ if it satisfies the
conditions $(1), (2)$ and \benum \item[(4)] $ \omega $ is the
natural antisymmetric bilinear form on $ \mathcal{A}_{1} \oplus
\mathcal{A}^{\ast}_{1} $ given by \beq \label{eq7} \omega(x +
a^{\ast}, y +  b^{\ast}) = -\langle x, b^{\ast} \rangle + \langle
a^{\ast}, y \rangle \mbox { for all } x, y \in \mathcal{A}_{1},
a^{\ast}, b^{\ast} \in \mathcal{A}^{\ast}_{1}, \eeq and $ \omega $
is a Connes cocycle on $ \mathcal{A} $. \eenum
\end{definition}

Let us now  give some notations useful in the sequel. Let $ \mathcal{A} $ be an associative algebra.

 Considering the representations of the  left $L$ and right $R$ multiplication operations defined as:
 \begin{eqnarray}
 L: \A & \longrightarrow & \mathfrak{gl}(\A)  \cr
  x  & \longmapsto & L_x:
  \begin{array}{ccc}
 \A &\longrightarrow & \A \cr 
  y & \longmapsto & x \cdot y, 
   \end{array}
\end{eqnarray}
\begin{eqnarray}
    R: \A & \longrightarrow & \mathfrak{gl}(\A)  \cr
     x  & \longmapsto & R_x:
     \begin{array}{ccc}
    \A &\longrightarrow & \A \cr 
     y & \longmapsto & y \cdot x,
      \end{array}
 \end{eqnarray}
 The dual maps $L^{*}, R^{*}$  of the linear maps $L, R,$ are defined, respectively, as: 
$\displaystyle L^{*}, R^{*}: \A \rightarrow \mathfrak{gl}(\A^{*})$
 such that:
 \newpage
\beq\label{dual1}
 L^*: \A & \longrightarrow & \mathfrak{gl}(\A^*)  \cr
  x  & \longmapsto & L^*_x:
      \begin{array}{llll}
 \A^* &\longrightarrow & \A^* \\ 
  u^* & \longmapsto & L^*_x u^*: 
      \begin{array}{llll}
\A  &\longrightarrow&  \K \cr
v  &\longmapsto& \left< L^{*}_xu^{*}, v \right>
 := \left<  u^{*}, L_x v\right>, 
      \end{array}
     \end{array}
 \eeq
\beq\label{dual2}
 R^*: \A & \longrightarrow & \mathfrak{gl}(\A^*)  \cr
  x  & \longmapsto & R^*_x:
      \begin{array}{llll}
 \A^* &\longrightarrow & \A^* \\ 
  u^* & \longmapsto & R^*_x u^*: 
      \begin{array}{lllll}
\A  &\longrightarrow&  \K \cr
v  &\longmapsto& \left< R^{*}_xu^{*}, v \right>
 := \left<  u^{*}, R_x v\right>, 
      \end{array}
     \end{array}
\eeq
for all $x, v \in \A, u^{*} \in \A^{*},$ where $\A^{*}$ is the dual space of $\A.$

 Let $ \sigma : \mathcal{A}\otimes \mathcal{A} \rightarrow \mathcal{A}\otimes \mathcal{A}  $ be the exchange operator defined as
\beqs
\sigma(x \otimes y) = y \otimes x,
\eeqs for all $ x, y \in \mathcal{A} $.

 An associative Yang-Baxter equation (AYBE) in the associative algebra $\mathcal{A}$ is defined by \cite{[C.Bai5]}
\beq \label{AYEB}
r_{12}r_{13} + r_{13}r_{23} - r_{23}r_{12} = 0,
\eeq
where $ r = \sum_{i} x_{i} \otimes y_{i} \in \mathcal{A} \otimes \mathcal{A}$ and
\beqs
r_{12}r_{13} = \sum_{i,j} x_{i}x_{j} \otimes y_{i} \otimes y_{j},
\eeqs
\beqs
r_{13}r_{23} = \sum_{i,j} x_{i} \otimes x_{j} \otimes y_{i}y_{j},
\eeqs
\beqs
r_{23}r_{12} = \sum_{i,j} x_{j} \otimes x_{i}y_{j} \otimes y_{i}.
\eeqs

\begin{definition}
Let $V_{1} $, $ V_{2} $ be two vector spaces. For a linear map $ \phi : V_{1} \rightarrow V_{2} $, we denote the dual (linear) map by $ \phi^{\ast} : V^{\ast}_{2} \rightarrow V^{\ast}_{1} $ given by 
\beqs
\langle v, \phi^{\ast}(u^{\ast})\rangle = \langle \phi(v), u^{\ast} \rangle
\eeqs 
for all  $ v \in V_{1} $, $ u^{\ast} \in V^{\ast}_{2} $.
\end{definition}

\begin{definition} 
Let $ \mathcal{A} $ be an associative algebra. An \textbf{antisymmetric infinitesimal
 bialgebra} structure on $ \mathcal{A} $ is a linear map $ \Delta: \mathcal{A} \rightarrow
\mathcal{A} \otimes \mathcal{A} $ such that
\begin{enumerate}
\item $  \Delta^{\ast} : \mathcal{A}^{\ast} \otimes \mathcal{A}^{\ast}  \rightarrow \mathcal{A}^{\ast}  $ defines an associative algebra structure on $ \mathcal{A}^{\ast} $;
\item $ \Delta $ satisfies the following equations:
\beq
\Delta (x\cdot y) = (id \otimes L(x))\Delta(y) + (R(y) \otimes id)\Delta(x),
\eeq
\beq
&&(L(y) \otimes id - id \otimes R(y))\Delta(x) \cr
&&+ \sigma [(L(x) \otimes id - id \otimes R(x))\Delta(y)] = 0,
\eeq
for all $ x, y \in \mathcal{A} $.
\end{enumerate}
\end{definition}
We denote this bialgebra structure by $ (\mathcal{A}, \Delta) $ or $(\mathcal{A},
 \mathcal{A}^{\ast})$.
\begin{theorem}  \label{c3}
Let $(\mathcal{A}, \cdot)$ and $(\mathcal{A}^{\ast}, \circ)$ be two associative algebras. Then, the following conditions are equivalent:
\begin{enumerate}
\item there is a double construction of a Frobenius algebra associated with  $(\mathcal{A}, \cdot) $ and $(\mathcal{A}^{\ast}, \circ)$;
\item $(\mathcal{A}, \mathcal{A}^{\ast})$ is an antisymmetric infinitesimal bialgebra.
\end{enumerate}
\end{theorem}

\begin{corollary}\label{c2}
Let $ \mathcal{A} $ be an associative algebra and $ r \in \mathcal{A} \otimes \mathcal{A} $. Suppose that $ r $ is antisymmetric. Then the map $ \Delta $ defined by
\begin{eqnarray}
\Delta(x) = (id \otimes L(x) - R(x) \otimes id)r \mbox { for all }  x \in \mathcal{A}
\end{eqnarray}
 induces an associative algebra structure on $ \mathcal{A}^{\ast} $
 such that $ (\mathcal{A}, \mathcal{A}^{\ast}) $ is an antisymmetric infinitesimal bialgebra if (\ref{AYEB}) is satisfied. 
\end{corollary}

\begin{proposition}  \label{pro2.4.4}
Let $(\mathcal{A},\cdot)$ be an associative algebra and let $ r
\in \mathcal{A} \otimes \mathcal{A} $ be an antisymmetric solution
of the associative Yang-Baxter equation in $ \mathcal{A} $. Then,
the corresponding
double construction of Frobenius algebra
$(\mathcal{AD}(\mathcal{A}),\ast) $ associated to $\mathcal{A}$
and $\mathcal A^*$
 is given from the
product in $ \mathcal{A} $ as follows: \beq \label{eq27} a^{\ast}
\ast b^{\ast} = a^{\ast} \circ b^{\ast} =
R^{\ast}(r(a^{\ast}))b^{\ast} + L^{\ast}(r(b^{\ast})) a^{\ast},
\eeq \beq \label{eq28} x \ast a^{\ast} = x \cdot r(a^{\ast}) -
r(R^{\ast}(x)a^{\ast}) + R^{\ast}(x)a^{\ast}, \eeq \beq
\label{eq29} a^{\ast} \ast x = r(a^{\ast}) \cdot x -
r(L^{\ast}(x)a^{\ast}) + L^{\ast}(x)a^{\ast}, \eeq for any $ x \in
\mathcal{A}, a^{\ast}, b^{\ast} \in \mathcal{A}^{\ast} $.
\end{proposition}

\begin{definition}
Let $ \mathcal{A} $ be a vector space
with two bilinear products denoted by $
\prec $ and $ \succ $. Then $ (\mathcal{A}, \prec, \succ) $ is
called a \textbf{dendriform algebra} if, for any $ x, y, z \in
\mathcal{A} $,
\begin{eqnarray*}
(x \prec y) \prec z &=& x \prec (y \ast z), \cr
(x \succ y) \prec z &=& x \succ (y \prec z), \cr
x \succ (y \succ z) &=& ( x \ast y) \succ z,
\end{eqnarray*}
where $ x \ast y = x \prec y + x \succ y $.
\end{definition}

Let $ (\mathcal{A}, \prec, \succ) $ be a dendriform algebra. For any $ x \in \mathcal{A} $, let $ L_{\succ}(x),  R_{\succ}(x) $ and $ L_{\prec}(x),$ \\$ R_{\prec}(x) $ denote the left and right multiplication operators of $(\mathcal{A}, \prec)$ and $(\mathcal{A}, \succ)$, respectively:
\beqs
L_{\succ}(x) y = x \succ y, R_{\succ}(x) y = y \succ x, L_{\prec}(x) y = x \prec y, 
R_{\prec}(x) y = y \prec x,
\eeqs
for all $ x, y \in \mathcal{A} $. Moreover, let $ L_{\succ}, R_{\succ}, L_{\prec}, R_{\prec} : \mathcal{A} \rightarrow gl(\mathcal{A}) $ be four linear maps with $ x \mapsto L_{\succ}(x),  x \mapsto R_{\succ}(x),  x \mapsto L_{\prec}(x), $ and $  x \mapsto R_{\prec}(x)  $, respectively. It is known that the product given by \cite{[N8]} 
\beq \label{eq32}
x \ast y = x \prec y + x \succ y, \mbox { for all } x, y \in \mathcal{A},
\eeq
defines an associative algebra. We call $ (\mathcal{A}, \ast) $ the associated associative algebra of $ (\mathcal{A}, \prec, \succ) $ and $ (\mathcal{A}, \succ, \prec) $ is called a compatible dendriform algebra structure on the associative algebra $ (\mathcal{A}, \ast) $.

\begin{theorem} \label{theo4.1.1}
Let $(\mathcal{A}, \ast)$ be an associative algebra and let $ \omega $ be a non-degenerate Connes cocycle. 
Then,  there exists a compatible dendriform algebra structure $ \succ, \prec $ on $  \mathcal{A} $ given by 
\beq \label{eq53}
\omega(x \succ y, z) = \omega(y, z \ast x), \ \ \ \omega(x \succ y, z) = \omega(x, y \ast z) \mbox { for all } x, y \in \mathcal{A}.
\eeq 
\end{theorem}

\begin{corollary} Let $(T(\mathcal{A}) = \mathcal{A} \bowtie \mathcal{A}^{\ast}, \omega)$ be a double construction of the Connes cocycle. Then, 
 there exists a compatible dendriform algebra structure $ \succ, \prec $ on $ T(\mathcal{A}) $ 
defined by the equation (\ref{eq53}).
  Moreover, $ \mathcal{A} $ and $ \mathcal{A}^{\ast} $, endowed with this product,
 are dendriform subalgebras. 
\end{corollary}

\begin{definition} 
Let $ \mathcal{A} $ be a vector space. A \textbf{dendriform D-bialgebra} structure on $ \mathcal{A} $ is a set of linear maps
$ (\Delta_{\prec}, \Delta_{\succ}, \beta_{\prec}, \beta_{\succ}) $  given by
 \ \ $ \Delta_{\prec}, \Delta_{\succ} : \mathcal{A} \rightarrow \mathcal{A} \otimes \mathcal{A} $, \ \
$ \beta_{\prec},  \beta_{\succ} : \mathcal{A}^{\ast} \rightarrow \mathcal{A}^{\ast} \otimes \mathcal{A}^{\ast} $, such that
\benum
\item[(a)] $ (\Delta^{\ast}_{\prec}, \Delta^{\ast}_{\succ}) : \mathcal{A}^{\ast} \otimes \mathcal{A}^{\ast} \rightarrow \mathcal{A}^{\ast} $
 defines a dendriform algebra structure $(\succ_{\mathcal{A}^{\ast}}, \prec_{\mathcal{A}^{\ast}}) $ on $ \mathcal{A}^{\ast} $;
\item[(b)] $ (\beta^{\ast}_{\prec}, \beta^{\ast}_{\succ}) : \mathcal{A} \otimes \mathcal{A} \rightarrow A $ defines a dendriform algebra
 structure $(\succ_{\mathcal{A}}, \prec_{\mathcal{A}})$ on $  \mathcal{A} $;
\item[(c)] the following equations are satisfied
\beq \label{eq55}
\Delta_{\prec}(x \ast_{\mathcal{A}} y) = (\id \otimes L_{\prec_{\mathcal{A}}}(x))\Delta_{\prec}(y) + (R_{\mathcal{A}}(y)\otimes \id)\Delta_{\prec}(y),
\eeq
\beq \label{eq56}
\Delta_{\succ}(x \ast_{\mathcal{A}} y) = (\id \otimes L_{\prec_{\mathcal{A}}}(x))\Delta_{\succ}(y) + (R_{\prec_{A}}(y)\otimes \id)\Delta_{\succ}(y),
\eeq
\beq \label{eq57}
\beta_{\prec}(a^{\ast} \ast_{\mathcal{A}^{\ast}} b^{\ast}) = (\id \otimes L_{\prec_{\mathcal{A}^{\ast}}}(a^{\ast}))\beta_{\prec}(b^{\ast}) + (R_{\mathcal{A}^{\ast}}(b^{\ast})\otimes \id)\beta_{\prec}(a^{\ast})
\eeq
\beq \label{eq58}
\beta_{\succ}(a^{\ast} \ast_{\mathcal{A}^{\ast}} b^{\ast}) = (\id \otimes L_{\mathcal{A}^{\ast}}(a^{\ast}))\beta_{\succ}(b^{\ast}) + (R_{\prec_{\mathcal{A}^{\ast}}}(b^{\ast})\otimes \id)\beta_{\succ}(a^{\ast}) ,
\eeq
\beq \label{eq59}
(L_{\mathcal{A}}(x)\otimes \id - \id \otimes R_{\prec_{\mathcal{A}}}(x))\Delta_{\prec}(y) + \sigma[(L_{\succ_{\mathcal{A}}}(y)\otimes (-\id)\otimes R_{\mathcal{A}}(y))\Delta_{\prec}(y)] = 0 ,
\eeq
\beq \label{eq60}
(L_{\mathcal{A}^{\ast}}(a^{\ast})\otimes \id - \id \otimes R_{\prec_{\mathcal{A}^{\ast}}}(a^{\ast}))\beta_{\prec}(b^{\ast}) 
+ \sigma[(L_{\succ_{\mathcal{A}^{\ast}}}(b^{\ast})\otimes(-id)\otimes R_{\mathcal{A}^{\ast}}(b^{\ast}))\beta_{\succ}(a^{\ast})] = 0,
\eeq
hold for any $ x, y \in \mathcal{A} $ and
$ a^{\ast}, b^{\ast} \in \mathcal{A}^{\ast} $,
where $ L_{\mathcal{A}} = L_{\succ_{\mathcal{A}}} + L_{\prec_{\mathcal{A}}}, R_{\mathcal{A}}
= R_{\succ_{\mathcal{A}}} + R_{\prec_{\mathcal{A}}}, L_{\mathcal{A}^{\ast}} = L_{\succ_{\mathcal{A}^{\ast}}}
+ L_{\prec_{\mathcal{A}^{\ast}}},  R_{\mathcal{A}^{\ast}} = R_{\succ_{\mathcal{A}^{\ast}}} + R_{\prec_{\mathcal{A}^{\ast}}} $.

\eenum
We also  denote it by $ (\mathcal{A}, \mathcal{A}^{\ast}, \Delta_{\succ}, \Delta_{\prec},  \beta_{\succ}, \beta_{\prec}) $
 or simply $(\mathcal{A}, \mathcal{A}^{\ast})$.
\end{definition}

\begin{theorem}
Let $(A, \prec_{A}, \succ_{A})$ and $(A^{\ast}, \prec_{A^{\ast}}, \succ_{A^{\ast}} )$ be two dendriform algebras. Let $(A, \ast_{A})$ and   $(A^{\ast}, \ast_{A^{\ast}})$ be the associated associative algebras, respectively. Then, the following
 conditions are equivalent:
 \begin{enumerate}
 \item there is a double construction of the Connes cocycle associated with
 $(A, \ast_{A})$ and $(A^{\ast}, \ast_{A^{\ast}})$;
\item $ (A, A^{\ast}) $ is a dendriform $ D- $bialgebra.
\end{enumerate}
\end{theorem}

\begin{corollary} \label{cor1}
Let $ (\mathcal{A}, \succ, \prec) $ be a dendriform algebra and
 $ r \in \mathcal{A} \otimes \mathcal{A} $. 
Suppose that
 $ r  $ is symmetric and $ r $ satisfies the equation
\begin{eqnarray}
\label{eqD}
r_{12} \ast r_{13} = r_{13} \prec r_{23} + r_{23} \succ r_{12}.
\end{eqnarray}
Then, the maps $  \Delta_{\succ}$ and $\Delta_{\prec} $ are  defined, respectively, by
\begin{eqnarray}
\Delta_{\succ}(x) &=& (id \otimes L(x) - R_{\prec}(x) \otimes id)r_{\succ}, \cr
\Delta_{\prec}(x) &=& (id \otimes L_{\succ}(x) - R(x) \otimes id)r_{\prec}, \forall x
\in \mathcal{A},
\end{eqnarray}
 where  $ r_{\succ} = - r$ and $ r_{\prec} = r$ induce a dendriform algebra structure on $ \mathcal{A}^{\ast} $ such that $(\mathcal{A}, \mathcal{A}^{\ast})$ is a dendriform $D$-bialgebra. Equation (\ref{eqD}) is called a \textbf{$D$-equation} in $ \mathcal{A} $.
\end{corollary}
\begin{proposition}  \label{pro4.4.6}
Let $ (\mathcal{A}, \succ, \prec) $ be a dendriform algebra  and let $ r \in \mathcal{A} \otimes \mathcal{A} $
 be a symmetric solution of the $D$-equation in $ \mathcal{A} $. Then,  the corresponding double
construction of Connes cocycle associated to $\mathcal{A}$ and
$\mathcal A^*$ is given from the products in $ \mathcal{A} $ as
follows:
 \beq \label{eq80} 
 a^{\ast} \prec b^{\ast} &=&
-R^{\ast}_{\succ}(r(a^{\ast}))b^{\ast} +
L^{\ast}(r(b^{\ast}))a^{\ast}, \cr 
a^{\ast} \succ b^{\ast} &=&
R^{\ast}(r(a^{\ast}))b^{\ast}
-L^{\ast}_{\prec}(r(b^{\ast}))a^{\ast}, \cr
a^{\ast} \ast
b^{\ast} &=& a^{\ast} \succ b^{\ast} + a^{\ast} \prec b^{\ast} =
R^{\ast}_{\prec}(r(a^{\ast}))b^{\ast} +
L^{\ast}_{\succ}(r(b^{\ast}))a^{\ast}, \cr
x \succ a^{\ast} &=& x
\succ r(a^{\ast}) - r(R^{\ast}(x)a^{\ast}) + R^{\ast}(x)a^{\ast},
\cr 
x \prec a^{\ast} &=&  x \prec r(a^{\ast}) +
r(R^{\ast}_{\succ}(x)a^{\ast}) - R^{\ast}_{\succ}(x)a^{\ast}, \cr
x \ast a^{\ast} &=& x \ast r(a^{\ast}) -
r(R^{\ast}_{\prec}(x)a^{\ast}) + R^{\ast}_{\prec}(x)a^{\ast}, \cr
a^{\ast} \succ x &=& r(a^{\ast}) \succ x +
r(L^{\ast}_{\prec}(x)a^{\ast}) -  L^{\ast}_{\prec}(x)a^{\ast}, \cr
a^{\ast} \prec x &=& r(a^{\ast}) \prec x -  r(L^{\ast}(x)a^{\ast})
+ L^{\ast}(x)a^{\ast}, \cr
a^{\ast} \ast x &=& r(a^{\ast}) \ast x
-  r(L^{\ast}_{\succ}(x)a^{\ast}) + L^{\ast}_{\succ}(x)a^{\ast}
\eeq for any $ x \in \mathcal{A}, a^{\ast}, b^{\ast} \in
\mathcal{A}^{\ast} $.
\end{proposition}

In the sequel, unless otherwise stated, all the parameters belong to the complex field $\mathbb{C}$.
 \section{$ 1 $-dimensional associative algebras}
 In this section, we investigate the solutions of the associative Yang-Baxter equation, dendriform algebras structures and classify the solutions of $D$-equations in the case of $ 1-$dimensional associative algebras.
\subsection{Solutions of the associative Yang-Baxter equation }
Let $(\mathcal{A}, \cdot)$ be an associative algebra with a basis $\lbrace e_{1} \rbrace$ and $ r= a_{11}e_{1}\otimes e_{1} \in \mathcal{A}\otimes \mathcal{A}.$ Then the AYBE becomes
\beqs
a^{2}_{11}(e_1\cdot  e_1\otimes e_1\otimes e_1 + e_1\otimes  e_1\otimes e_1\cdot e_1 - e_1\otimes e_1\cdot e_1\otimes e_1)=0.
\eeqs
\begin{proposition}
There are only two non-isomorphic $1-$dimensional associative algebras. The solutions of the corresponding associative Yang-Baxter equation  are given in  Table 1.

\begin{center}
{\bf Table 1}:
 Solutions of the $ 1- $dimensional  associative Yang-Baxter equation.

\begin{tabular}{c  c}
\hline
Associative algebra $ \mathcal{A} $  & Solutions of the AYBE
  \\ \hline
 $ \mathcal{A}_{1}: e_{1}\cdot e_{1}= 0 $ & $ r= a_{11}e_{1}\otimes e_{1} $
\\ \hline
$ \mathcal{A}_{2}: e_{1}\cdot e_{1}= e_{1}$ & $ r= 0 $
\\ \hline
\end{tabular}
\end{center}
\end{proposition}

\subsubsection{Antisymmetric solutions and Frobenius algebra structures}
Using the Proposition \ref{pro2.4.4}, we obtain the results presented in Table $ 2  $.
\begin{center}
{\bf Table 2}: Antisymmetric solutions and Frobenius algebra structures of  $ 1- $dimensional  associative algebras.

\begin{tabular}{c c c}
\hline
Associative algebra $ \mathcal{A}$  & Antisymmetric solutions  & Frobenius algebra   structures over $ \mathcal{A}\oplus  \mathcal{A}^{\ast} $
  \\ \hline
  $ \mathcal{A}_{1} $  & $  r=0  $ & $ e_{1}\ast e_{1}= e_{1}\ast e^{\ast}_{1}= e^{\ast}_{1}\ast e_{1}= 0$
  \\ \hline
$ \mathcal{A}_{2} $  & $ r=0  $ & $ e_{1}\ast e_{1}= e_{1};e_{1}\ast e^{\ast}_{1}= e^{\ast}_{1} $;
    $e^{\ast}_{1}\ast e_{1}= e^{\ast}_{1}$

\\ \hline
\end{tabular}
\end{center}

\subsection{Dendriform algebra structures and  classification of solutions of  the
$D$-equation}
\begin{proposition} The compatible dendriform algebra structures  in $ 1 $-dimensional  associative algebras and related solutions of D-equations  are given in  Table 3.

\begin{center}
{\bf Table 3}: The $ 1- $dimensional dendriform algebras and
classification of solutions of $D$-equations.

\begin{tabular}{c  c  c c}
\hline
Associative  & Dendriform  & Solutions of
 \cr
 algebra  $ \mathcal{A} $   & algebra structures    &   $D$-equation
    \\ \hline
$ \mathcal{A}_{1} $: $ e_{1}\cdot e_{1}= 0 $ & $ D^{1}_{1} $: $ e_{1}\prec e_{1}= e_{1}\succ e_{1}=0  $ & $ r= a_{11}e_{1}\otimes e_{1} $
\\ \hline
  $ \mathcal{A}_{2} $: $ e_{1}\cdot e_{1}= e_{1} $ & $ D^{2}_{1} $: $ e_{1}\succ e_{1}= \lambda e_{1} $, $ e_{1}\prec e_{1}= (1-\lambda)e_{1}$;   $ \lambda =0, 1 $  &
    $ r= a_{11}e_{1}\otimes e_{1} $
 \\ \hline
\end{tabular}
\end{center}
\end{proposition}
\textbf{Proof }
Let us consider the associative algebra $ \mathcal{A}_{2} $. We  set
$e_{1}\succ e_{1}= ae_{1}; e_{1}\prec e_{1}= be_{1} $. Since $ e_{1}\cdot e_{1}= e_{1} $, then $ a + b= 1 $. Moreover, we have the equations
\beqs
&&(e_{1}\prec e_{1})\prec e_{1}= e_{1}\prec(e_{1}\cdot e_{1}), \cr
&&(e_{1}\succ e_{1})\prec e_{1}= e_{1}\succ(e_{1}\prec e_{1}), \cr
&& e_{1}\succ(e_{1}\succ e_{1})= (e_{1}\cdot e_{1})\succ e_{1}
\eeqs
which give  $ ab= 0; ab= ba $ and $ a(a- 1)= 0 $, respectively. Then, we obtain $ a=0$ and $ b= 1 $ or $ a= 1 $ and $ b= 0 $ yielding the compatible dendriform algebra structures  on $ \mathcal{A}_{2} $.
The solutions of D-equations in these dendriform algebras are given by direct computation.
$ \hfill \square $

\subsubsection{Symmetric solutions and Connes cocycles structures}
Using the Proposition \ref{pro4.4.6}, we obtain the following results in Table $ 4 $, where
$ \displaystyle D^j_i,$  $i, j\in \mathbb{N}^{\ast}$,  means the $i-th$ dendriform class associated with the $j-th$ class of associative algebra.

\begin{center}
{\bf Table 4}:
 Symmetric solutions and Connes cocycles  of  $ 1- $dimensional  associative algebras.

\begin{tabular}{c c c}
\hline
Dendriform  algebra  & Symmetric solutions &  Connes cocycles over $ \mathcal{A}\oplus  \mathcal{A}^{\ast} $
\\ \hline
$ D^{1}_{1} $ & $ a_{11}e_{1}\otimes e_{1}  $ & $ e_{1}\ast e_{1}=e_{1}\ast e^{\ast}_{1}= e^{\ast}_{1}\ast e_{1}=e^{\ast}_{1}\ast e^{\ast}_{1}= 0$
\\ \hline
$ D^{2}_{1} $   & $ a_{11}e_{1}\otimes e_{1}  $  & ,$ e^{\ast}_{1}\ast e^{\ast}_{1}= -a_{11} e^{\ast}_{1}$, $  e^{\ast}_{1}\ast e_{1}= a_{11}(\lambda- 1)e_{1} +\lambda e^{\ast}_{1}  $ \cr
        & &$ e_{1}\ast e_{1}= e_{1}$, $  e_{1}\ast e^{\ast}_{1}= -a_{11}\lambda e_{1} + (1 -\lambda)e^{\ast}_{1};   \lambda= 0, 1 $
\\ \hline
\end{tabular}
\end{center}
 \section{$ 2 $-dimensional associative algebras}

 Using a similar approach as in the previous section, we now consider $  2-$dimensional associative
  algebras. Let $(\mathcal{A}, \cdot)$ be an associative algebra with a basis $\lbrace e_{1}, e_{2}
   \rbrace$ and $ r = \sum^{2}_{i, j=1} a_{ij}e_{i}\otimes e_{j} \in \mathcal{A} \otimes \mathcal{A}.$ Then, the AYBE (\ref{AYEB}), i.e. $r_{12}r_{13} + r_{13}r_{23} - r_{23}r_{12} = 0,$ is satisfied for
\beqs
r_{12}= \sum^{2}_{i, j=1} a_{ij}e_{i}\otimes e_{j}\otimes 1,
\eeqs
\beqs
r_{13}= \sum^{2}_{i, j=1} a_{ij}e_{i}\otimes 1\otimes e_{j},
\eeqs
\beqs
r_{23}= \sum^{2}_{i, j=1} a_{ij} 1\otimes e_{i}\otimes e_{j}.
\eeqs

\subsection{Solutions of the associative Yang-Baxter equation}
 The classification of 2-dimensional complex pre-Lie algebras, including the classification of 2-dimensional complex associative algebras, was performed in \cite{[Bu]}. Then, the $2$-dimensional complex associative algebras can be split into  $7$ classes \cite{[C.Bai6]}.
\begin{proposition} \label{pro3}
The solutions of the associative Yang-Baxter equation (\ref{AYEB})
 in 2-dimensional associative algebras are given
in  Table 5.

\begin{center}
{\bf Table 5}:
 Solutions of the $ 2- $dimensional  associative Yang-Baxter equation.
\begin{tabular}{c c}
\hline
Associative algebra $ \mathcal{A} $ & Solutions of the AYBE
 \\ \hline
 $ \mathcal{A}_{1}:$ $ e_{1}\cdot e_{1}= e_{2} $  &
 $ \left(
\begin{array}{cc}
0 & 0  \\
0 & a_{22}  \\
\end{array}
\right)$
\\ \hline
$ \mathcal{A}_{2}:e_{1}\cdot e_{1}= e_{1}, e_{1}\cdot e_{2}= e_{2}  $ & $ \left(
\begin{array}{cc}
0 & a_{12}  \\
0 & a_{22}  \\
\end{array}
\right) $; $ \left(
\begin{array}{cc}
0 & a_{12}  \\
a_{21} & 0  \\
\end{array}
\right) a_{12} = - a_{21}\neq 0   $
\\ \hline
$ \mathcal{A}_{3}:e_{1}\cdot e_{1}= e_{1}, e_{2}\cdot e_{1}= e_{2}  $ & $ \left(
\begin{array}{cc}
0 & 0  \\
a_{21} & a_{22}  \\
\end{array}
\right) $; $ \left(
\begin{array}{cc}
0 & a_{12}  \\
a_{21} & 0  \\
\end{array}
\right) a_{12} = - a_{21}\neq 0   $
\\ \hline
$ \mathcal{A}_{4}:e_{1}\cdot e_{1}= e_{1}, e_{1}\cdot e_{2}= e_{2}  $, & $ \left(
\begin{array}{cc}
0 & a_{12}  \\
0 & 0  \\
\end{array}
\right) $; $ \left(
\begin{array}{cc}
0 & 0  \\
a_{21} & 0  \\
\end{array}
\right) a_{21}\neq 0   $; \cr
$ e_{2}\cdot e_{1}= e_{2} $ & $ \left(
\begin{array}{cc}
0 & 0  \\
0 & a_{22}  \\
\end{array}
\right) a_{22}\neq 0   $
\\ \hline
 $ \mathcal{A}_{5}:e_{i}\cdot e_{j}= 0; i, j=1, 2$ & $ \left(
\begin{array}{cc}
a_{11} & a_{12}  \\
a_{21} & a_{22}  \\
\end{array}
\right)$
\\ \hline
 $ \mathcal{A}_{6}:$ $ e_{2}\cdot e_{2}= e_{2} $ &
  $ \left(
\begin{array}{cc}
a_{11} & 0  \\
0 & 0  \\
\end{array}
\right)$
\\ \hline
 $ \mathcal{A}_{7}:e_{1}\cdot e_{1}= e_{1} $; $ e_{2}\cdot e_{2}= e_{2} $ & $ \left(
\begin{array}{cc}
0 & 0  \\
0 & 0  \\
\end{array}
\right)$
\\ \hline
\end{tabular}
\end{center}
\end{proposition}
\textbf{Proof }
We set $ r = \sum_{i, j}a_{ij}e_{i}\otimes e_{j} \in \mathcal{A}\otimes \mathcal{A}$, where $ i, j=1, 2 $ and $ a_{ij} \in \mathbb{C} $. Then, by direct computation of (\ref{AYEB}),
we get the results.
$ \hfill \square $

\subsubsection{Antisymmetric solutions and Frobenius algebra structures}
Using the Proposition \ref{pro2.4.4}, we get the  results of  Table $ 6 $.

\begin{center}
{\bf Table 6}:
 Antisymmetric solutions and Frobenius algebra structures of $ 2-$dimensional  associative algebras.

\begin{tabular}{c c c}
\hline
Associative  & Antisymmetric  & Frobenius algebra \cr
  algebra $ \mathcal{A}$  & solutions     &  structures over $ \mathcal{A}\oplus  \mathcal{A}^{\ast} $
  \\ \hline
  $ \mathcal{A}_{1} $  & $ r =0  $ & $ e_{1}\ast e_{1}=  e_{2};e_{1}\ast e^{\ast}_{2}= e^{\ast}_{2} $; $e^{\ast}_{2}\ast e_{1}= e^{\ast}_{1}$
    \\ \hline
$ \mathcal{A}_{2} $  & $ a_{12}e_{1}\otimes e_{2}  $ & $ e_{1}\ast e_{1}=  e_{1};e_{1}\ast e_{2}= e_{2} $; $e_{1}\ast e^{\ast}_{1}= e^{\ast}_{1};$ $  e^{\ast}_{2}\ast e_{1}= e^{\ast}_{2}  $ \cr
  & $ + a_{21}e_{2}\otimes e_{1} $;  &  $e^{\ast}_{2}\ast e^{\ast}_{2}= -a_{12} e^{\ast}_{2}$;  $e^{\ast}_{2}\ast e^{\ast}_{1}= a_{21}e^{\ast}_{1}$; $ e_{1}\ast e^{\ast}_{2}= a_{21}e_{1}$  \cr
                         &  $ a_{12} = - a_{21} \neq 0 $  &  $ e_{2}\ast e^{\ast}_{2}= -a_{12}e_{2} + e^{\ast}_{1} $; $ e^{\ast}_{2}\ast e_{2}= -a_{12}e_{2} $; $ e^{\ast}_{1}\ast e_{1}= a_{21}e_{2} + e^{\ast}_{1} $
\\ \hline
   $ \mathcal{A}_{3} $  & $ a_{12}e_{1}\otimes e_{2}  $ & $ e_{1}\ast e_{1}=  e_{1};e_{2}\ast e_{1}= e_{2} $; $e_{1}\ast e^{\ast}_{2}= e^{\ast}_{2}$;  $  e^{\ast}_{1}\ast e_{1}= e^{\ast}_{1}  $ \cr
  & $ + a_{21}e_{2}\otimes e_{1} $;  & $ e^{\ast}_{2}\ast e^{\ast}_{2}= a_{21} e^{\ast}_{2}$;  $e^{\ast}_{1}\ast e^{\ast}_{2}= a_{21}e^{\ast}_{1} $; $ e_{1}\ast e^{\ast}_{1}= -a_{12}e_{2} + e^{\ast}_{1} $ \cr
 &  $ a_{12} = - a_{21} \neq 0 $  & $ e_{2}\ast e^{\ast}_{2}= a_{21}e_{2}$;
 $ e^{\ast}_{2}\ast e_{1}= -a_{12}e_{1} $; $ e^{\ast}_{2}\ast e_{2}= a_{21}e_{2} + e^{\ast}_{1} $
  \\ \hline
    $ \mathcal{A}_{4} $  &  & $ e_{1}\ast e_{1}=  e_{1};e_{1}\ast e_{2}= e_{2} $; $e_{2}\ast e_{1}= e_{2};  e^{\ast}_{1}\ast e_{1}=  e^{\ast}_{1}$\cr
    & $ r =0  $ & $ e_{1}\ast e^{\ast}_{1}= e^{\ast}_{1} $
     $e^{\ast}_{2}\ast e_{1}= e^{\ast}_{2} $; $e^{\ast}_{2}\ast e_{2}= e^{\ast}_{1}$
    \\ \hline
$ \mathcal{A}_{5} $  & $ a_{12}e_{1}\otimes e_{2}  $ & $ e_{i}\ast e_{j}= e^{\ast}_{i}\ast e^{\ast}_{j}=  e_{i}\ast e^{\ast}_{j}= 0$ \cr
  & $ + a_{21}e_{2}\otimes e_{1} $;  &   $e^{\ast}_{i}\ast e_{j}= 0$; $i, j= 1,2$  \cr
                         &  $ a_{12} = - a_{21} \neq 0 $  &
  \\ \hline
    $ \mathcal{A}_{6} $  & $ r =0  $ & $ e_{2}\ast e_{2}=  e_{2};e_{2}\ast e^{\ast}_{2}= e^{\ast}_{2} $; $e^{\ast}_{2}\ast e_{2}= e^{\ast}_{2}$
    \\ \hline
   $ \mathcal{A}_{7} $  &  & $ e_{1}\ast e_{1}= e_{1}$; $e_{2}\ast e_{2}= e_{2}$; $ e_{2}\ast e^{\ast}_{2}= e^{\ast}_{2}$ \cr
  &  $ r =0  $  &  $e^{\ast}_{1}\ast e_{1}= e^{\ast}_{1}$; $ e^{\ast}_{2}\ast e_{2}= e^{\ast}_{2}$;  $e_{1}\ast e^{\ast}_{1}= e^{\ast}_{1}$
  \\ \hline
\end{tabular}
\end{center}
\subsection{Dendriform algebra structures   and  classification of solutions of  the
$D$-equation} The classification of $2-$dimensional complex
dendriform algebras was firstly studied in \cite{[C.Bai7]}, but
unfortunately with some mistakes. In fact, there exists a natural
anti-isomorphism between dendriform algebras \beqs F(x\succ_{1}
y)=F(y)\prec_{2} F(x), F(x\prec_{1} y)=F(y)\succ_{2} F(x). \eeqs
So, the dendriform algebras appear in terms of pairs. For example, consider the following dendriform algebras given by:
\beqs
&& D_{1}: e_{1}\prec_{1} e_{1}= e_{1}, e_{1}\prec_{1} e_{2}= e_{2}; \cr
&& D_{2}: e_{1}\succ_{2} e_{1}= e_{1}, e_{2}\succ_{2} e_{1}= e_{2}; 
\eeqs
and the map $ F $ defined by$ F(e_{1})= e_{1}, F(e_{2})= e_{2}.$ We have 
$F(e_{1}\prec_{1} e_{1})=e_{1}=F(e_{1})\succ_{2} F(e_{1}); F(e_{1}\prec_{1} e_{2})=e_{2}=F(e_{2})\succ_{2} F(e_{1}).$
Therefore, there exists an anti-isomrphism between the dendriform algebras $ D_{1} $ and $ D_{2}.$
 By this property, one can  classify dendriform algebras. For that, let $  e_{i}, e_{j} \in \mathcal{A}; i, j=1, 2$. Then,
\beq \label{Q}
e_{i}\cdot e_{j}= e_{i}\prec e_{j} + e_{i}\succ e_{j}, \mbox { where }
 e_{i}\succ e_{j} = \sum^{2}_{k= 1} a^{k}_{ij}e_{k},
e_{i}\prec e_{j} = \sum^{2}_{k= 1} b^{k}_{ij}e_{k}. \eeq Computing
the equation (\ref{Q}) with the condition $ e_{i}\cdot
e_{j}=\sum^{2}_{k= 1} (a^{k}_{ij} +b^{k}_{ij} ) e_{k}$ , we get
the  compatible dendriform algebra structures  on the associative
algebra. The solutions of $D$-equations in these dendriform
algebras are then determined by direct computation. There results  the following: 
 
\begin{proposition}
 The compatible dendriform algebra structures  in $2$-dimensional associative algebras
and the solutions of related  D-equations   are given in  Table 7.
\begin{center}
{\bf Table 7}:
 $ 2- $dimensional dendriform algebras and  classification of solutions of the
$D$-equation.
\begin{tabular}{c c c c}
\hline
Associative  & Dendriform  & Solutions of
 \cr
 algebra   & algebra     &   $D$-equation \cr
 $ \mathcal{A} $  & structures  &
    \\ \hline
$ \mathcal{A}_{1} $: & $ D^{1}_{1, \lambda} $: $ e_{1}\prec e_{1}=\lambda e_{2}  $ & $ \left(
\begin{array}{cc}
0 & 0  \\
0 & a_{22}  \\
\end{array}
\right),
 $ \cr
$ e_{1}\cdot e_{1}= e_{2} $  &  $  e_{1}\succ e_{1}=(1- \lambda) e_{2} $  & $ \left(
\begin{array}{cc}
0 & 0  \\
a_{21} & a_{22}  \\
\end{array}
\right) \lambda= 0, a_{21}\neq 0
$ \cr
  &   & $ \left(
\begin{array}{cc}
0 & a_{12}  \\
a_{21} & a_{22}  \\
\end{array}
\right) a_{12}\neq a_{21}
$ \cr
  &  & $ \lambda= \frac{a^{2}_{12}- a_{12}a_{21}}{a^{2}_{21}- a_{12}a_{21}} $
\\ \hline
  $ \mathcal{A}_{2} $: & $ D^{2}_{1} $: $ e_{1}\succ e_{2}= e_{2} $ & $  \left(
\begin{array}{cc}
a_{11} & 0  \\
a_{21} & 0  \\
\end{array}
\right);  \left(
\begin{array}{cc}
0 & 0  \\
0 & a_{22}  \\
\end{array}
\right), a_{22}\neq 0$ \cr
$ e_{1}\cdot e_{1}= e_{1} $  &  $  e_{1}\succ e_{1}=  e_{1} $  & $ \left(
\begin{array}{cc}
0 & 0  \\
a_{21} & a_{22}  \\
\end{array}
\right) a_{21}, a_{22}\neq 0$ \cr
 $ e_{1}\cdot e_{2}= e_{2} $ &  & $ \left(
\begin{array}{cc}
a_{11} & a_{12}  \\
a_{21} & a_{22}  \\
\end{array}
\right) a_{12}\neq 0$ \cr
  &  & $ a_{21}= \frac{a_{11}a_{22}}{a_{12}} $ \cr

\cline{2-3}  & $ D^{2}_{2} $: $  e_{1}\succ e_{2}= e_{2} $  &  $ \left(
\begin{array}{cc}
0 & a_{12}  \\
0 & a_{22}  \\
\end{array}
\right)  $ \cr
& $  e_{1}\prec e_{1}=  e_{1} $  &  $ \left(
\begin{array}{cc}
a_{11} & 0  \\
0 & 0  \\
\end{array}
\right) $  \cr
                         &   & $ \left(
\begin{array}{cc}
0 & a_{12}  \\
a_{21} & a_{22}  \\
\end{array}
\right)a_{21}= a_{12}\neq 0 $ \cr
\cline{2-3}
     & $ D^{2}_{3} $: $ e_{1}\prec e_{1}= e_{1} $ & $ \left(
\begin{array}{cc}
0 & 0  \\
0 & a_{22}  \\
\end{array}
\right) $ \cr
 & $ e_{1}\prec e_{2}= e_{2}  $ & $ \left(
\begin{array}{cc}
a_{11} & a_{12}  \\
a_{21} & a_{22}  \\
\end{array}
\right) a_{12}= a_{21}, a_{11}\neq 0 $  \cr
& & $ a_{22}=\frac{a^{2}_{12}}{a_{11}} $
 \\ \hline
\end{tabular}
\end{center}

\begin{center}
\begin{tabular}{c c c c}
\hline
Associative  & Dendriform  & Solutions of
 \cr
 algebra   & algebra     &   $D$-equation \cr
 $ \mathcal{A} $  & structures  &
    \\ \hline
 & $ D^{2}_{4} $: $ e_{2}\prec e_{1}= e_{2} $ & $ \left(
\begin{array}{cc}
a_{11} & 0  \\
a_{21} & 0  \\
\end{array}
\right) $ \cr
 & $ e_{1}\succ e_{2}= e_{2}  $ & $ \left(
\begin{array}{cc}
a_{11} & a_{12}  \\
0 & 0  \\
\end{array}
\right) a_{12}\neq 0 $ \cr
 & $ e_{1}\prec e_{1}= e_{1} $ & $ \left(
\begin{array}{cc}
a_{11} & a_{12}  \\
a_{21} & a_{22}  \\
\end{array}
\right), a_{21}= a_{11}; a_{12}= a_{22}\neq 0 $  \cr
 &$ e_{2}\succ e_{1}= -e_{2} $ & 
    \\ \hline
& $ D^{2}_{5} $: $ e_{1}\succ e_{2}= e_{2} $ & $ \left(
\begin{array}{cc}
0 & a_{12}  \\
0 & a_{22}  \\
\end{array}
\right)  $ \cr
& $ e_{1}\succ e_{1}= - e_{2} $ & $ \left(
\begin{array}{cc}
a_{11} & a_{12}  \\
a_{21} & a_{22}  \\
\end{array}
\right) a_{12}= a_{21}\neq 0$ \cr
& $ e_{1}\prec e_{1}=  e_{1} +  e_{2} $ &
    \\ \hline
          $ \mathcal{A}_{3} $:& $ D^{3}_{1} $: $  e_{1}\succ e_{1}= e_{1} $  & $\left(
\begin{array}{cc}
a_{11} & 0  \\
a_{21} & 0  \\
\end{array}
\right) $ \cr
& $ e_{2}\prec e_{1}= e_{2} $ & $ \left(
\begin{array}{cc}
0 & a_{12}  \\
a_{21} & 0  \\
\end{array}
\right) a_{21}= a_{12}\neq 0
$  \cr
&& $ \left(
\begin{array}{cc}
0 & a_{12}  \\
a_{21} & a_{22}  \\
\end{array}
\right) a_{21}= a_{12}, a_{22}\neq 0
$ \cr
\cline{2-3}
 &$ D^{3}_{2} $: $  e_{1}\prec e_{1}= e_{1} $& $ \left(
\begin{array}{cc}
0 & 0  \\
0 & a_{22}  \\
\end{array}
\right)  $ \cr
& $ e_{2}\prec e_{1}= e_{2} $ & $ \left(
\begin{array}{cc}
a_{11} & a_{12}  \\
a_{21} & a_{22}  \\
\end{array}
\right) a_{21}= a_{12}; a_{22}=\dfrac{a^{2}_{12}}{a_{11}}, a_{11}\neq 0  $ \cr
    \cline{2-3}
  & $ D^{3}_{3} $: $ e_{1}\succ e_{1}= e_{1} - e_{2} $   & $ \left(
\begin{array}{cc}
0 & 0  \\
0 & a_{22}  \\
\end{array}
\right);  \left(
\begin{array}{cc}
0 & 0  \\
a_{21} & a_{22}  \\
\end{array}
\right)
$ \cr
  & $ e_{1}\prec e_{1}= e_{2} $  & $ a_{22}= -a_{21}\neq 0 $ \cr
  & $  e_{2}\prec e_{1}= e_{2} $ & $ \left(
\begin{array}{cc}
0 & a_{12}  \\
a_{21} & a_{22}  \\
\end{array}
\right),
 a_{12}=a_{21} \neq 0 $ \cr
  &  & $ \left(
\begin{array}{cc}
a_{11} & a_{12}  \\
a_{21} & a_{22}  \\
\end{array}
\right) a_{22}= - a_{21}
,$\cr
& & $ a_{12}= -a_{11} $ \cr
\cline{2-3}
    & $ D^{3}_{4} $: $ e_{1}\succ e_{2}=  e_{2}$ & $ \left(
\begin{array}{cc}
0 & a_{12}  \\
0 & a_{22}  \\
\end{array}
\right) $ \cr
& $ e_{1}\succ e_{1}= e_{1}$ & $ \left(
\begin{array}{cc}
0 & 0  \\
a_{21} & a_{22}  \\
\end{array}
\right) $ \cr
& $ e_{1}\prec e_{2}= -e_{2}$ & $ \left(
\begin{array}{cc}
a_{11} & a_{12}  \\
a_{21} & a_{22}  \\
\end{array}
\right) a_{11}\neq 0; a_{22}= \dfrac{a_{21}a_{12}}{a_{11}} $  \cr
& $ e_{2}\prec e_{1}=  e_{2} $ &  
    \\ \hline
 $ \mathcal{A}_{4} $: & $ D^{4}_{1} $: $ e_{1}\succ e_{2}=  e_{2}  $ & $  \left(
\begin{array}{cc}
a_{11} & 0  \\
0 & a_{22}  \\
\end{array}
\right);
 $ \cr
$ e_{1}\cdot e_{1}= e_{1} $ & $ e_{1}\prec e_{1}= e_{1}  $ & $ \left(
\begin{array}{cc}
0 & a_{12}  \\
0 & 0  \\
\end{array}
\right)  a_{12}\neq 0;$ \cr
 $ e_{1}\cdot e_{2}= e_{2} $ & $ e_{2}\prec e_{1}= e_{2} $   & $ \left(
\begin{array}{cc}
a_{11} & a_{12}  \\
0 & 0  \\
\end{array}
\right) $ \cr
$ e_{2}\cdot e_{1}= e_{2} $ & & $ a_{11}, a_{12} \neq 0 $\cr
\cline{2-3}
  & $ D^{4}_{2} $:  $  e_{1}\succ e_{2}= e_{2} $  &  $ \left(
\begin{array}{cc}
0 & 0  \\
0 & a_{22}  \\
\end{array}
\right)$ \cr
&  $ e_{1}\succ e_{1}= e_{1} $ & $ \left(
\begin{array}{cc}
0 & a_{12}  \\
0 & 0  \\
\end{array}
\right), a_{12}\neq 0$  \cr
                         & $  e_{2}\succ e_{1}=  e_{2} $   & $ \left(
\begin{array}{cc}
0 & a_{12}  \\
a_{21} & a_{22}  \\
\end{array}
\right)a_{21}= a_{12}\neq 0 $ \cr
& & $ \left(
\begin{array}{cc}
a_{11} & a_{12}  \\
0 & 0  \\
\end{array}
\right)a_{11} \neq 0 $
  \\ \hline
\end{tabular}
\end{center}

\begin{tabular}{c c c c}
\hline
Associative  & Dendriform  & Solutions of
 \cr
 algebra   & algebra     &   $D$-equation \cr
 $ \mathcal{A} $  & structures  &
    \\ \hline
 & $ D^{4}_{3} $: $ e_{1}\succ e_{2}= e_{2} $ & $ \left(
\begin{array}{cc}
a_{11} & 0  \\
0 & a_{22}  \\
\end{array}
\right)$ \cr
& $ e_{1}\succ e_{1}= e_{1} $ & $ \left(
\begin{array}{cc}
0 & 0  \\
a_{21} & 0  \\
\end{array}
\right), a_{21}\neq 0  $ \cr
& $ e_{2} \prec e_{1}= e_{2} $ & $ \left(
\begin{array}{cc}
0 & a_{12}  \\
0 & 0  \\
\end{array}
\right) a_{12} \neq 0   $ \cr
    \cline{2-3}
 & $ D^{4}_{4} $: $ e_{1}\succ e_{2}=  e_{2}  $ & $  \left(
\begin{array}{cc}
a_{11} & 0  \\
0 & a_{22}  \\
\end{array}
\right);
 $ \cr
& $ e_{1}\succ e_{1}= -e_{2}  $ & $ \left(
\begin{array}{cc}
a_{11} & a_{12}  \\
0 & 0  \\
\end{array}
\right) a_{12}\neq 0 $ \cr
  & $ e_{2}\succ e_{1}=  e_{2} $   &  \cr
 & $ e_{1}\prec e_{1}= e_{1} + e_{2}  $ & 
 \\ \hline
      $ \mathcal{A}_{5} $:  & $ D^{5}_{1} $: $ e_{2}\succ e_{2}= \beta e_{1} $ & $ \left(
\begin{array}{cc}
a_{11} & a_{12}  \\
a_{21} & 0  \\
\end{array}
\right) \beta=0;$ $  \left(
\begin{array}{cc}
a_{11} & a_{12}  \\
0 & 0  \\
\end{array}
\right) \beta=1$ \cr
  $ e_{i}\cdot e_{j}= 0;  $ &  $  e_{2}\prec e_{2}= -\beta e_{1} $ & $  \left(
\begin{array}{cc}
a_{11} & 0  \\
a_{21} & 0  \\
\end{array}
\right) a_{21}\neq 0, \beta=1$ \cr
$ i, j= 1, 2 $ & $ \beta=0, 1 $ &  $  \left(
\begin{array}{cc}
a_{11} & a_{12}  \\
a_{21} & a_{22}  \\
\end{array}
\right) a_{22}\neq 0; \beta= 0 $
  \\ \hline
$ \mathcal{A}_{6} $:$ e_{2}\cdot e_{2}= e_{2} $ & $ D^{6}_{1} $: $ e_{2}\succ e_{2}= e_{2} $ & $ \left(
\begin{array}{cc}
a_{11} & a_{12}  \\
0 & 0  \\
\end{array}
\right) $;
 $  \left(
\begin{array}{cc}
a_{11} & 0  \\
0 & a_{22}  \\
\end{array}
\right) a_{22}\neq 0 $ \cr
\cline{2-3}
 & $ D^{6}_{2} $: $ e_{2}\succ e_{2}=  e_{2}  $ & $ \left(
\begin{array}{cc}
0 & a_{12} \\
a_{21} & 0  \\
\end{array}
\right) $ \cr
& $ e_{2}\succ e_{1}= e_{1}$ & $ \left(
\begin{array}{cc}
0 & 0 \\
0 & a_{22}  \\
\end{array}
\right) a_{22}\neq 0 $ \cr
  & $ e_{2}\prec e_{1}= -e_{1} $   & $ \left(
\begin{array}{cc}
a_{11} & a_{12} \\
a_{21} & 0  \\
\end{array}
\right) a_{11}\neq 0; a_{21}= a_{12} $  \cr
  \cline{2-3}
   & $ D^{6}_{3} $: $ e_{2}\prec e_{2}=  e_{2}  $ & $ \left(
\begin{array}{cc}
a_{11} & a_{12} \\
0 & 0  \\
\end{array}
\right) $ \cr
& $ e_{1}\prec e_{2}= e_{1}   $ & $ \left(
\begin{array}{cc}
a_{11} & a_{12} \\
a_{21} & 0  \\
\end{array}
\right) a_{21}= a_{12}\neq 0 $ \cr
  & $ e_{1}\succ e_{2}= -e_{1} $   & $ \left(
\begin{array}{cc}
0 & 0 \\
0 & a_{22}  \\
\end{array}
\right) a_{22}\neq 0 $ \cr
 \cline{2-3}
    & $ D^{6}_{4} $: $ e_{1}\succ e_{1}=  e_{2} $ & $ \left(
\begin{array}{cc}
0 & a_{12}  \\
0 & 0  \\
\end{array}
\right) $ \cr
& $ e_{2}\succ e_{2}=  e_{2} $  & \cr
& $ e_{1}\prec e_{1}= -e_{2} $  & $ \left(
\begin{array}{cc}
0 & 0  \\
0 & a_{22}  \\
\end{array}
\right), a_{22}\neq 0  $
\\ \hline

 $ \mathcal{A}_{7} $:  & $ D^{7}_{1} $: $ e_{1}\prec e_{2}=  e_{1}  $ & $ \left(
\begin{array}{cc}
a_{11} & 0 \\
0 & 0  \\
\end{array}
\right) $; $ \left(
\begin{array}{cc}
a_{11} & 0 \\
a_{21} & 0  \\
\end{array}
\right)a_{21}\neq 0 $   \cr
$ e_{1}\cdot e_{1}= e_{1}  $ & $ e_{2}\prec e_{2}= e_{2}  $ & $ \left(
\begin{array}{cc}
0 & 0 \\
0 & a_{22}  \\
\end{array}
\right) a_{22}\neq 0 $ \cr
  & $ e_{1}\succ e_{1}= e_{1} $   & $ \left(
\begin{array}{cc}
a_{11} & a_{12} \\
a_{21} & a_{22}  \\
\end{array}
\right), a_{21}\neq a_{12}$ \cr
$ e_{2}\cdot e_{2}= e_{2}  $ &  $ e_{1}\succ e_{2}= -e_{1}  $ &  $ \left(
\begin{array}{cc}
a_{11} & a_{12} \\
a_{21} & a_{22}  \\
\end{array}
\right)$ $a_{21}=a_{12}\neq 0, a_{12}= a_{22}$ \cr
&& $ \left(
\begin{array}{cc}
0 & a_{12}  \\
0 & a_{22} \\
\end{array}
\right) a_{12}= a_{22} \neq 0$ 
\\ \hline
\end{tabular}

\begin{center}
\begin{tabular}{c c c c}
\hline
Associative  & Dendriform  & Solutions of
 \cr
 algebra   & algebra    &   $D$-equation \cr
 $ \mathcal{A} $  & structures  &
    \\ \hline
  & $ D^{7}_{2} $: $ e_{1}\succ e_{1}= e_{1} $ & $ \left(
\begin{array}{cc}
0 & 0  \\
a_{21} & a_{22}  \\
\end{array}
\right);$ $ \left(
\begin{array}{cc}
0 & a_{12}  \\
a_{21} & a_{22}  \\
\end{array}
\right)a_{12}= a_{22}\neq 0 $  \cr
& $ e_{2}\succ e_{2}= e_{2} $ & $ \left(
\begin{array}{cc}
a_{11} & 0  \\
0 & a_{22}  \\
\end{array}
\right)a_{11}\neq 0$  \cr
 && $ \left(
\begin{array}{cc}
a_{11} & a_{12}  \\
a_{21} & a_{22}  \\
\end{array}
\right)a_{11}= a_{21}\neq 0 $ $a_{22}= a_{12}$\cr
\cline{2-3}
  & $ D^{7}_{3} $: $ e_{2}\prec e_{2}= e_{2} $ & $ \left(
\begin{array}{cc}
a_{11} & 0  \\
0 & a_{22}  \\
\end{array}
\right) $ \cr
& $ e_{1}\succ e_{1}= e_{1} $ & $ \left(
\begin{array}{cc}
a_{11} & 0  \\
a_{21} & 0  \\
\end{array}
\right), a_{21}= a_{11}\neq 0 $ \cr 
\cline{2-3}
 & $ D^{7}_{4} $: $ e_{1}\prec e_{2}= -e_{2}  $ & $ \left(
\begin{array}{cc}
0 & 0 \\
0 & a_{22}  \\
\end{array}
\right);$ $\left(
\begin{array}{cc}
a_{11} & 0 \\
0 & 0  \\
\end{array}
\right) a_{11}\neq 0 $  \cr
& $ e_{2}\prec e_{2}= e_{2}  $ & $ \left(
\begin{array}{cc}
a_{11} & a_{12} \\
0 & 0  \\
\end{array}
\right) a_{11}=a_{12}\neq 0$ \cr
  & $ e_{1}\succ e_{2}= e_{2} $ &$\left(
\begin{array}{cc}
a_{11} & a_{12} \\
a_{21} & 0  \\
\end{array}
\right) a_{21}= a_{11}\neq 0,$ $ a_{12}\neq 0$  \cr
 & $ e_{1}\succ e_{1}= e_{1}  $ & $\left(
\begin{array}{cc}
a_{11} & 0  \\
a_{21} & 0 \\
\end{array}
\right) a_{21}= a_{11}\neq 0$ 
    \\ \hline
    & $ D^{7}_{5, \lambda} $: $ e_{1}\prec e_{1}=  e_{1} $ & $ \left(
\begin{array}{cc}
a_{11} & 0  \\
0 & a_{22}  \\
\end{array}
\right)\lambda =0;$ $ \left(
\begin{array}{cc}
a_{11} & 0  \\
0 & 0  \\
\end{array}
\right)\lambda\neq 0 $ \cr
& $ e_{2}\prec e_{1}= -\lambda e_{2} $  &$ \left(
\begin{array}{cc}
0 & 0  \\
0 & a_{22}  \\
\end{array}
\right)\lambda , a_{22}\neq 0$\cr
& $ e_{2}\succ e_{1}= \lambda e_{2}$  & $ \left(
\begin{array}{cc}
a_{11} & 0  \\
a_{21} & a_{22}  \\
\end{array}
\right), \lambda = -1,$  $ a_{21}= a_{11}\neq 0$ \cr
& $ e_{2}\succ e_{2}= e_{2} $& $ \left(
\begin{array}{cc}
a_{11} & a_{12}  \\
0 & a_{22}  \\
\end{array}
\right)a_{21}= a_{11}\neq 0,$ $ \lambda= -1, a_{12}\neq 0 $ \cr
&& $ \left(
\begin{array}{cc}
0 & a_{12}  \\
0 & a_{22}  \\
\end{array}
\right) a_{22}= a_{12}\neq 0 $\cr
  \cline{2-3}
 & $ D^{7}_{6} $: $ e_{1}\prec e_{1}= e_{1}  $ & $ \left(
\begin{array}{cc}
0 & 0 \\
0 & a_{22}  \\
\end{array}
\right);$ $ \left(
\begin{array}{cc}
0 & 0 \\
a_{21} & a_{22}  \\
\end{array}
\right)a_{22}= a_{21}\neq 0 $ \cr
& $ e_{2}\prec e_{1}= - e_{1} $ & $ \left(
\begin{array}{cc}
0 & a_{12} \\
0 & a_{22}  \\
\end{array}
\right)a_{22}= a_{12}\neq 0$\cr
  & $ e_{2}\succ e_{1}=  e_{1} $   &  $ \left(
\begin{array}{cc}
a_{11} & a_{12} \\
a_{21} & a_{22}  \\
\end{array}
\right) a_{11}\neq 0;$ $ a_{22}= a_{21}= a_{12} $\cr
  & $ e_{2}\succ e_{2}=  e_{2} $ & $ \left(
\begin{array}{cc}
0 & a_{12} \\
a_{21} & a_{22}  \\
\end{array}
\right)$ $ a_{22}= a_{12}= a_{21}\neq 0 $ \cr
\cline{2-3}
  & $ D^{7}_{7} $: $ e_{1}\prec e_{1}= e_{1} $ & $ \left(
\begin{array}{cc}
a_{11} & 0  \\
0 & 0  \\
\end{array}
\right) $ \cr
& $ e_{2}\prec e_{1}= -e_{1} $ &  \cr
& $ e_{2}\prec e_{2}= e_{1} + e_{2} $ & \cr
& $ e_{2}\succ e_{1}= e_{1} $ & \cr
& $ e_{2}\succ e_{2}= -e_{1} $ &\cr
\cline{2-3}
 & $ D^{7}_{8, \lambda} $: $ e_{1}\prec e_{1}= \lambda e_{2}  $ & $ \left(
\begin{array}{cc}
0 & 0 \\
0 & a_{22}  \\
\end{array}
\right) $  \cr
& $ e_{1}\succ e_{1}= e_{1} -\lambda e_{2}  $ & $ \left(
\begin{array}{cc}
0 & a_{12} \\
0 & a_{22}  \\
\end{array}
\right) a_{12}= a_{22}\neq 0 $ \cr
  & $ e_{2}\succ e_{2}= e_{2} $   &  \cr
  & $ \lambda \neq 0 $ &  
\\ \hline
\end{tabular}
\end{center}

\begin{center}
\begin{tabular}{c c c c}
\hline
Associative  & Dendriform  & Solutions of
 \cr
 algebra   & algebra    &   $D$-equation \cr
 $ \mathcal{A} $  & structures  &
\\ \hline
   & $ D^{7}_{9, \lambda} $: $ e_{1}\prec e_{1}= \lambda e_{2}  $ & $ \left(
\begin{array}{cc}
0 & 0 \\
0 & a_{22}  \\
\end{array}
\right) $  \cr
& $ e_{2}\prec e_{2}= e_{2}  $ &  \cr
  & $ e_{1}\succ e_{1}= e_{1}- \lambda e_{2} $   &  \cr
  & $ \lambda \neq 0 $ & \cr
\cline{2-3}
& $ D^{7}_{10, \lambda} $: $ e_{1}\prec e_{2}= - e_{2}  $ & $ \left(
\begin{array}{cc}
0 & 0 \\
0 & a_{22}  \\
\end{array}
\right) $  \cr
& $ e_{1}\prec e_{1}= e_{1} -\lambda e_{2}  $ & $ \left(
\begin{array}{cc}
0 & a_{12} \\
0 & 0  \\
\end{array}
\right) a_{12}\neq 0 $ \cr
  & $ e_{2}\prec e_{1}= -e_{2} $   &   \cr
& $ e_{1}\succ e_{2}= e_{2} $  & \cr
&$ e_{2}\succ e_{1}= e_{2} $ & \cr
& $ e_{1}\succ e_{1}= \lambda e_{2} $ & \cr
& $ e_{2}\succ e_{2}= e_{2} $ & \cr
& $ \lambda\neq 0  $ & \cr
\cline{2-3}
 & $ D^{7}_{11, \lambda} $: $ e_{1}\prec e_{2}= - e_{2}  $ & $ \left(
\begin{array}{cc}
0 & 0 \\
0 & a_{22}  \\
\end{array}
\right) $  \cr
& $ e_{1}\prec e_{1}= e_{1} -\lambda e_{2}  $ & $ \left(
\begin{array}{cc}
0 & a_{12} \\
0 & 0  \\
\end{array}
\right) a_{12}\neq 0 $ \cr
  & $ e_{2}\prec e_{1}= -e_{2} $   &   \cr
& $ e_{2}\prec e_{2}= e_{2} $  & \cr
&$ e_{1}\succ e_{2}= e_{2} $ & \cr
& $ e_{2}\succ e_{1}=  e_{2} $ & \cr
& $ e_{1}\succ e_{1}= \lambda e_{2} $ & \cr
& $ \lambda\neq 0 $ &
\\ \hline
\end{tabular}
\end{center}
\end{proposition}

\subsubsection{Symmetric solutions and Connes cocycles structures}
Using the Proposition \ref{pro4.4.6}, we obtain the results in Table $ 8 $, where
$ \displaystyle D^j_i $, with $i, j\in \mathbb{N}^{\ast}$, means the $i-th$ dendriform class associated with the $j-th$ class of associative algebra.

{\bf Table 8}:
 Symmetric solutions and Connes cocycles structures of  $ 2- $dimensional  associative algebras.

\begin{tabular}{c c c}
\hline
Dendriform   & Symmetric  &  Connes cocycles \cr
 algebra   & solutions     &  structures over $ \mathcal{A}\oplus  \mathcal{A}^{\ast} $\cr
structures &&
    \\ \hline
$ D^{1}_{1, \lambda} $ & $ a_{22}e_{2}\otimes e_{2}  $ & $ e_{1}\ast e_{1}=  e_{2}$; $e_{1}\ast e^{\ast}_{2}= \lambda e^{\ast}_{1}$; $e^{\ast}_{2}\ast e_{1}= (1- \lambda)e^{\ast}_{1}$
\\ \hline
$ D^{2}_{1} $   & $ a_{11}e_{1}\otimes e_{1}  $  & $ e_{1}\ast e_{1}= e_{1}$; $ e_{1}\ast e_{2}= e_{2};$ $ e^{\ast}_{2}\ast e_{1}= e^{\ast}_{2};$ $e^{\ast}_{1}\ast e_{2}= -a_{11} e_{2}$\cr
                      &    &  $  e_{1}\ast e^{\ast}_{1}= -a_{11}e_{1}$;
                         $ e^{\ast}_{2}\ast e^{\ast}_{1}= -a_{11} e^{\ast}_{2} $;
                        $  e^{\ast}_{1}\ast e_{1}= e^{\ast}_{1}   $\cr
                        & &
                        $  e^{\ast}_{1}\ast e^{\ast}_{1}= -a_{11} e^{\ast}_{1} $ \cr
\cline{2-3}
    & $ a_{22}e_{2}\otimes e_{2}  $  & $ e_{1}\ast e_{1}= e_{1}$; $ e_{1}\ast e_{2}= e_{2};$ $ e^{\ast}_{2}\ast e_{1}= e^{\ast}_{2} $; \cr
                         &   $ a_{22}\neq 0 $ &
                         $  e_{1}\ast e^{\ast}_{2}= -a_{22}e_{2};$  
                        $ e^{\ast}_{1}\ast e_{1}= e^{\ast}_{1}  $ \cr
 \cline{2-3} & $ a_{11}e_{1}\otimes e_{1} +  $  & $ e_{1}\ast e_{1}=  e_{1}; e_{1}\ast e_{2}= e_{2}; e^{\ast}_{1}\ast e^{\ast}_{1}= -a_{11}e^{\ast}_{1} $  \cr
 & $ a_{12}e_{1}\otimes e_{2} + $   &  $ e_{1}\ast e^{\ast}_{2}= -a_{12}e_{1}$; $ e^{\ast}_{2}\ast e^{\ast}_{2}=  a_{22}e^{\ast}_{1};$ $ e_{1}\ast e^{\ast}_{1}= -a_{11}e_{1} -a_{12}e_{2}  $ \cr
 &    $ a_{12}e_{2}\otimes e_{1} + $    & $ e^{\ast}_{1}\ast e_{1}=  a_{12}e_{2} + e^{\ast}_{1};$ $ e^{\ast}_{2}\ast e_{1}=  a_{22}e_{2} + e^{\ast}_{2};$  \cr
& $ a_{22}e_{2}\otimes e_{2}  $ & $ a^{2}_{12}= a_{11}a_{22} $
\\ \hline
   $ D^{2}_{2} $ & $ a_{12}e_{1}\otimes e_{2} + $ & $ e_{1}\ast e_{1}= e_{1}; e_{1}\ast e_{2}= e_{2};$ $ e^{\ast}_{1}\ast e^{\ast}_{2}= -a_{12}e^{\ast}_{1};$ $e^{\ast}_{2}\ast e^{\ast}_{1}= -a_{12}e^{\ast}_{1} $ \cr
  & $  a_{12}e_{2}\otimes e_{1} + $  & $ e^{\ast}_{2}\ast e^{\ast}_{2}= -a_{22}e^{\ast}_{1}$; $ e^{\ast}_{2}\ast e_{2}= -a_{12}e_{2}$; $ e_{1}\ast e^{\ast}_{1}= -a_{12}e_{2} + e^{\ast}_{1} $   \cr
                         & $ a_{22}e_{2}\otimes e_{2} $   & $   e_{1}\ast e^{\ast}_{2}= -a_{22}e_{2} -a_{12} e_{1};$ $  e^{\ast}_{2}\ast e_{1}= a_{22}e_{2} + e^{\ast}_{2};$  \cr
 &$ a_{12}\neq 0 $& 
\\ \hline
\end{tabular}

\begin{tabular}{c c c}
\hline
Dendriform   & Symmetric  &  Connes cocycles \cr
 algebra   & solutions     &  structures over $ \mathcal{A}\oplus  \mathcal{A}^{\ast} $\cr
structures &&
    \\ \hline
$ D^{2}_{2} $ & $ a_{11}e_{1}\otimes e_{1}  $  & $ e_{1}\ast e_{1}= e_{1}; e_{1}\ast e_{2}=e_{2}$; $ e^{\ast}_{2}\ast e_{1}= e^{\ast}_{2}; e^{\ast}_{1}\ast e^{\ast}_{1}= -a_{11}e^{\ast}_{1} $ \cr
 &  &  $ e^{\ast}_{1}\ast e_{1}= -a_{11}e^{\ast}_{1}$;  $  e_{1}\ast e^{\ast}_{1}= e^{\ast}_{1}$ \cr
 \cline{2-3} & $ a_{22}e_{2}\otimes e_{2}  $  & $ e_{1}\ast e_{1}= e_{1}; e_{1}\ast e_{2}=e_{2}; e_{1}\ast e^{\ast}_{1}= e^{\ast}_{1};$ $e_{1}\ast e^{\ast}_{2}= -a_{22}e_{2} $ \cr
 &  & $  e^{\ast}_{2}\ast e_{1}= a_{22}e_{2} +  e^{\ast}_{2}  $    
    \\ \hline
$ D^{2}_{3} $ & $ a_{11}e_{1}\otimes e_{1} + $ & $ e_{1}\ast e_{1}= e_{1}; e_{1}\ast e_{2}=e_{2};$ $ e^{\ast}_{2}\ast e^{\ast}_{1}= a_{12}e^{\ast}_{1}$  \cr
 & $ a_{12}e_{1}\otimes e_{2} + $ & $ e^{\ast}_{1}\ast e^{\ast}_{2}= -a_{12}e^{\ast}_{1}; e_{1}\ast e^{\ast}_{1}=e^{\ast}_{1}; e^{\ast}_{1}\ast e^{\ast}_{1}= -a_{11}e^{\ast}_{1}$ \cr
 & $  a_{21}e_{2}\otimes e_{1}    $ & $e^{\ast}_{1}\ast e_{1}= -a_{11}e_{1};$ $ e^{\ast}_{1}\ast e_{2}= -a_{11}e_{2}, e^{\ast}_{2}\ast e_{1}= -a_{12}e_{1} $ \cr
 & $ a_{22}e_{2}\otimes e_{2}  $ & $ e^{\ast}_{2}\ast e_{2}= -a_{12}e_{2};$
  $ e_{1}\ast e^{\ast}_{2}= -a_{22}e_{2} -a_{12}e_{1}$; \cr
 & $ a_{21}= a_{12} $ &   $ e_{2}\ast e^{\ast}_{2}= a_{12}e_{2} +a_{11}e_{1} + e^{\ast}_{1};  e^{\ast}_{2}\ast e^{\ast}_{2}= -a_{22}e^{\ast}_{1}$
\\ \hline    
$ D^{2}_{4} $ & $a_{11}e_{1}\otimes e_{1}$ &  $ e_{1}\ast e_{1}=  e_{1}$; $ e_{1}\ast e_{2}= e_{2}$; $ e_{1}\ast e^{\ast}_{1}= e^{\ast}_{1};  e_{1}\ast e^{\ast}_{2}= e^{\ast}_{2} $ \cr
&& $  e^{\ast}_{1}\ast e_{1}= -a_{11}e_{1} $; $e^{\ast}_{1}\ast e_{2}= -a_{11}e_{2} $; $  e^{\ast}_{2}\ast e_{1}= e^{\ast}_{2}; e^{\ast}_{1}\ast e^{\ast}_{1}= - a_{11}e^{\ast}_{1}$\cr
&& $  e^{\ast}_{1}\ast e^{\ast}_{2}= - a_{11}e^{\ast}_{2};  e^{\ast}_{2}\ast e^{\ast}_{1}= - a_{11}e^{\ast}_{2};  e^{\ast}_{2}\ast e_{2}= - a_{11}e_{1} -e^{\ast}_{1}$
\\ \hline
 $ D^{2}_{5} $  & $ a_{11}e_{1}\otimes e_{1} + $ & $ e_{1}\ast e_{1}= e_{1}; e_{1}\ast e_{2}=e_{2};$ $ e^{\ast}_{2}\ast e^{\ast}_{2}= -a_{12}e^{\ast}_{2}$ \cr
 & $ a_{12}e_{1}\otimes e_{2} + $ &$e^{\ast}_{1}\ast e^{\ast}_{1}= -a_{11}e^{\ast}_{1};$ $ e_{1}\ast e^{\ast}_{1}= e^{\ast}_{1}$;  $ e^{\ast}_{1}\ast e_{1}= -a_{11}e_{1} $ \cr
   &  $ a_{12}e_{2}\otimes e_{1}  $ & $ e^{\ast}_{1}\ast e_{2}= -a_{11}e_{2}$;  $ e^{\ast}_{2}\ast e_{2}= -a_{12}e_{2} $ \cr
   & $ a_{22}e_{2}\otimes e_{2}$ & $ e^{\ast}_{2}\ast e_{1}= (-a_{12} + a_{22})e_{2} -$
   $ e^{\ast}_{1} + e^{\ast}_{2} $ \cr
   & $ a_{21}=a_{12}\neq 0  $ & $ e_{1}\ast e^{\ast}_{2}= (a_{12} - a_{22})e_{2} +   e^{\ast}_{1} + (a_{11} - a_{12})e_{1} $ \cr
   &   & $ e^{\ast}_{2}\ast e^{\ast}_{1}= (a_{11} - a_{12})e^{\ast}_{1} -a_{11}e^{\ast}_{2} $ \cr
   \cline{2-3}
   & $ a_{22}e_{2}\otimes e_{2}  $ & $ e_{1}\ast e_{1}= e_{1}; e_{1}\ast e_{2}= e_{2}$; $e_{1}\ast e^{\ast}_{1}=  e^{\ast}_{1}$ \cr
 &&  $e^{\ast}_{1}\ast e_{1}= e^{\ast}_{1};$ $ e_{1}\ast e^{\ast}_{2}= -a_{22}e_{2} +  e^{\ast}_{1}$  \cr
 &&  $e^{\ast}_{2}\ast e_{1}= a_{22}e_{2} + e^{\ast}_{2} - e^{\ast}_{1} $
\\ \hline
 $ D^{3}_{1} $ & $ a_{22}e_{2}\otimes e_{2}  $ & $ e_{1}\ast e_{1}=  e_{1}$; $ e_{2}\ast e_{1}= e_{2}$; $ e^{\ast}_{1}\ast e_{1}= e^{\ast}_{1};  e^{\ast}_{2}\ast e_{1}= -a_{22}e_{2} $ \cr
& $ a_{22}\neq 0 $ &  $  e_{1}\ast e^{\ast}_{2}= e^{\ast}_{2} + a_{22}e_{2}$\cr
\cline{2-3}
& $ a_{11}e_{1}\otimes e_{1}  $ & $ e_{1}\ast e_{1}=  e_{1}$; $ e_{2}\ast e_{1}= e_{2}$; $ e^{\ast}_{1}\ast e_{1}= e^{\ast}_{1};  e_{1}\ast e^{\ast}_{2}= e^{\ast}_{2} $ \cr
&  &  $  e_{1}\ast e^{\ast}_{1}= -a_{11}e_{1}; e_{2}\ast e^{\ast}_{1}= -a_{11}e_{2}$ \cr
&  &  $  e^{\ast}_{1}\ast e^{\ast}_{1}= -a_{11}e^{\ast}_{1}; e^{\ast}_{1}\ast e^{\ast}_{2}= -a_{11}e^{\ast}_{2}$ \cr
\cline{2-3}
& $ a_{12}e_{1}\otimes e_{2} + $ &  \cr
& $ a_{21}e_{2}\otimes e_{1}  $ & \cr
& $ a_{12}= a_{21}\neq 0$ &
\\ \hline
 $ D^{3}_{2} $ & $ a_{22}e_{2}\otimes e_{2}$ & $ e_{1}\ast e_{1}=  e_{1}$; $ e_{2}\ast e_{1}= e_{2}$; $ e_{1}\ast e^{\ast}_{1}= e^{\ast}_{1}$ \cr
&  &  $ e^{\ast}_{2}\ast e_{1}= -a_{22}e_{2}; e_{1}\ast e^{\ast}_{2}= e^{\ast}_{2} + a_{22}e_{2}$\cr
\cline{2-3}
& $ a_{11}e_{1}\otimes e_{1}  $ & $ e_{1}\ast e_{1}=  e_{1}$; $ e_{2}\ast e_{1}= e_{2}$; $ e_{1}\ast e^{\ast}_{1}= e^{\ast}_{1};  e_{1}\ast e^{\ast}_{2}= e^{\ast}_{2} $ \cr
&  &  $  e_{2}\ast e^{\ast}_{1}= -a_{11}e_{2}; e^{\ast}_{1}\ast e_{1}= -a_{11}e_{1}$ 
\\ \hline
 $ D^{3}_{3} $ & $ a_{12}e_{1}\otimes e_{2} + $ & $ e_{1}\ast e_{1}=  e_{1}$; $ e_{2}\ast e_{1}= e_{2};$ $ e_{1}\ast e^{\ast}_{1}= -a_{11}e_{1}$; $ e_{2}\ast e^{\ast}_{1}= -a_{11}e_{2}$ \cr
&  $ a_{12}e_{2}\otimes e_{1} +  $  & $ e_{2}\ast e^{\ast}_{2}= -a_{12}e_{2}$; $ e^{\ast}_{1}\ast e^{\ast}_{1}= -a_{11}e^{\ast}_{1};$ $e^{\ast}_{2}\ast e^{\ast}_{1}=  a_{11}e^{\ast}_{1}$ \cr
&   $ a_{11}e_{1}\otimes e_{1} +$  & $ e^{\ast}_{2}\ast e^{\ast}_{2}= -a_{12}e^{\ast}_{2};$ $ e_{1}\ast e^{\ast}_{2}= a_{11}e_{1} + (a_{12} + a_{22})e_{2}  $ + $  e^{\ast}_{1} +e^{\ast}_{2}$    \cr
& $ a_{22}e_{2}\otimes e_{2} $ &  $ e^{\ast}_{1}\ast e_{1}=   a_{11}e_{2} + e^{\ast}_{1};$  $e^{\ast}_{1}\ast e^{\ast}_{2}= (-a_{11}  - a_{12})e^{\ast}_{1} $  $- a_{11} e^{\ast}_{2}  $\cr
& $ a_{22}=  a_{11}  $ & $e^{\ast}_{2}\ast e_{1}=(-a_{11}- a_{12})e_{1} + $ $ (- a_{12}- a_{22})e_{2}  - e^{\ast}_{1} $ \cr
\cline{2-3}
     $ D^{3}_{3} $ & $ a_{22}e_{2}\otimes e_{2}  $ & $ e_{1}\ast e_{1}=  e_{1}$; $ e_{2}\ast e_{1}= e_{2}$; $ e_{1}\ast e^{\ast}_{1}= e^{\ast}_{1}  $; $ e_{2}\ast e^{\ast}_{2}= -e^{\ast}_{1}  $   \cr
&& $ e^{\ast}_{2}\ast e^{\ast}_{2}= a_{22}e^{\ast}_{1}  $ ; $  e^{\ast}_{2}\ast e_{1}= - e^{\ast}_{1} + e^{\ast}_{2}   $ ; $  e_{1}\ast e^{\ast}_{2}= a_{22}e_{2} + $ $ e^{\ast}_{1}  + e^{\ast}_{2}  $\cr
\cline{2-3}
     $ D^{3}_{4} $ & $ a_{22}e_{2}\otimes e_{2}  $ & $ e_{1}\ast e_{1}=  e_{1}$; $ e_{2}\ast e_{1}= e_{2}$; $ e^{\ast}_{2}\ast e^{\ast}_{2}= a_{22}e^{\ast}_{1}$; $ e^{\ast}_{2}\ast e_{1}= e^{\ast}_{2}  $   \cr
&& $ e^{\ast}_{1}\ast e_{1}= e^{\ast}_{1}$ ; $e_{2}\ast e^{\ast}_{2}= - e^{\ast}_{1}$ ; $e_{1}\ast e^{\ast}_{2}= a_{22}e_{2} + e^{\ast}_{2}$ \cr
\cline{2-3}
& $ a_{11}e_{1}\otimes e_{1} + $  & $ e_{1}\ast e_{1}=  e_{1}$; $ e_{2}\ast e_{1}= e_{2}$; $ e^{\ast}_{1}\ast e^{\ast}_{1}= -a_{11}e^{\ast}_{1}$; $ e^{\ast}_{2}\ast e^{\ast}_{1}= -a_{11}e^{\ast}_{2}  $   \cr
& $ a_{12}e_{1}\otimes e_{2} + $  & $ e^{\ast}_{2}\ast e^{\ast}_{2}= (a_{22} -a_{12})e^{\ast}_{1} -2a_{12}e^{\ast}_{2}; e^{\ast}_{1}\ast e^{\ast}_{2}= -a_{11}e^{\ast}_{2}$ \cr
& $ a_{12}e_{2}\otimes e_{1} +  $  & $ e_{2}\ast e^{\ast}_{1}= -a_{11}e_{2}; e_{1}\ast e^{\ast}_{2}= -a_{12}e_{1} + (a_{12} + a_{22})e_{2} + e^{\ast}_{2} $ \cr
& $ a_{22}e_{2}\otimes e_{2}  $  & $  e_{2}\ast e^{\ast}_{2}= -a_{11}e_{1} -2a_{12}e_{2}-e^{\ast}_{1}; e_{1}\ast e^{\ast}_{1}= -a_{11}e_{1}  $ \cr
& $ a_{11}a_{22}= a^{2}_{12}; a_{11}\neq 0 $  & $ e^{\ast}_{2}\ast e_{1}= -a_{12}e_{2} + e^{\ast}_{2}; e^{\ast}_{1}\ast e_{1}= a_{12}e_{1} + e^{\ast}_{1} $
\\ \hline
$ D^{4}_{1} $ & $ a_{11}e_{1}\otimes e_{1} +  $ & $ e_{1}\ast e_{1}=  e_{1}$; $ e_{1}\ast e_{2}= e_{2};$ $ e_{2}\ast e_{1}= e_{2}; e^{\ast}_{2}\ast e_{1}= e^{\ast}_{2}$ \cr
  &  $ a_{22}e_{2}\otimes e_{2} $ & $ e^{\ast}_{1}\ast e^{\ast}_{1}= -a_{11} e^{\ast}_{1}$; $ e^{\ast}_{1}\ast e^{\ast}_{2}= -a_{11}e^{\ast}_{2}; e^{\ast}_{1}\ast e^{\ast}_{2}= -a_{11}e^{\ast}_{2}$   \cr
  &&$ e_{1}\ast e^{\ast}_{2}= e^{\ast}_{2};$ $ e^{\ast}_{2}\ast e^{\ast}_{1}= -a_{11} e^{\ast}_{2}$; $ e_{2}\ast e^{\ast}_{1}= -a_{11}e_{2}$ \cr
  &&  $ e_{1}\ast e^{\ast}_{1}= e^{\ast}_{1}$; $ e^{\ast}_{1}\ast e_{1}= -a_{11}e_{1}$
  \\ \hline
\end{tabular}

\begin{tabular}{c c c}
\hline
Dendriform   & Symmetric  &  Connes cocycles \cr
 algebra   & solutions     &  structures over $ \mathcal{A}\oplus  \mathcal{A}^{\ast} $ \cr
structures &&
 \\ \hline
 $ D^{4}_{2} $ & $ a_{11}e_{1}\otimes e_{1}   $ & $ e_{1}\ast e_{1}=  e_{1}$; $ e_{1}\ast e_{2}= e_{2}$;  $ e^{\ast}_{2}\ast e_{2}= a_{11} e_{1} + e^{\ast}_{1}$ \cr
  &&     $ e^{\ast}_{2}\ast e_{1}= e^{\ast}_{2}; e_{1}\ast e^{\ast}_{1}= -a_{11}e_{1} $ \cr
  && $ e^{\ast}_{1}\ast e^{\ast}_{1}= -a_{11} e^{\ast}_{1}$; $ e^{\ast}_{2}\ast e^{\ast}_{1}= -a_{11}e^{\ast}_{2}; e_{2}\ast e_{1}= e_{2}$\cr
  &   & $ e_{2}\ast e^{\ast}_{1}= -a_{11} e_{2}$; $ e^{\ast}_{1}\ast e_{2}= -a_{11}e_{2}; e^{\ast}_{1}\ast e_{1}= e^{\ast}_{1}$ \cr
\cline{2-3}
 & $ a_{22}e_{2}\otimes e_{2} + $ & $ e_{1}\ast e_{1}=  e_{1}$;$ e_{1}\ast e^{\ast}_{1}= -a_{12} e_{2}; e_{2}\ast e^{\ast}_{2}= -a_{12} e_{2}$  \cr
  & $ a_{12}e_{1}\otimes e_{2} + $  &  $ e_{2}\ast e_{1}= e_{2}; e^{\ast}_{2}\ast e_{2}= e^{\ast}_{1};$ $ e_{1}\ast e^{\ast}_{2}=  -a_{12}e_{1} - a_{22}e_{2} $ \cr
  & $ a_{21}e_{2}\otimes e_{1}  $ & $e^{\ast}_{2}\ast e_{1}=  e^{\ast}_{2} + a_{12}e_{2};$  $ e^{\ast}_{1}\ast e_{1}=  e^{\ast}_{1} - a_{12}e_{2} $  \cr
  &  &$ e_{1}\ast e_{2}= e_{2};$ $ e^{\ast}_{1}\ast e^{\ast}_{2}= -a_{12} e^{\ast}_{1}; e^{\ast}_{2}\ast e^{\ast}_{1}= -2a_{12} e^{\ast}_{1}$\cr
&& $ e^{\ast}_{2}\ast e^{\ast}_{2}=  (a_{12} - a_{22})e^{\ast}_{1} - a_{12}e^{\ast}_{2} $
\\ \hline
 $ D^{4}_{3} $ & $ a_{11}e_{1}\otimes e_{1}  + $ & $ e_{1}\ast e_{1}=  e_{1}$; $ e_{1}\ast e_{2}= e_{2};$ $e^{\ast}_{1}\ast e^{\ast}_{1}= -a_{11}e^{\ast}_{1}$; $ e^{\ast}_{1}\ast e^{\ast}_{2}= -a_{11}e^{\ast}_{2}$ \cr
  & $ a_{22}e_{2}\otimes e_{2}$  &  $ e_{2}\ast e_{1}= e_{2}; e^{\ast}_{1}\ast e_{1}= e^{\ast}_{1};$ $e^{\ast}_{2}\ast e^{\ast}_{1}= -a_{11}e^{\ast}_{2}$; $ e_{1}\ast e^{\ast}_{1}= -a_{11}e_{1}$ \cr
  &&     $ e_{1}\ast e^{\ast}_{2}= e^{\ast}_{2}; e^{\ast}_{2}\ast e_{1}= e^{\ast}_{2};$ $e_{2}\ast e^{\ast}_{1}= -a_{11}e^{\ast}_{2}$; $ e^{\ast}_{1}\ast e_{2}= -a_{11}e_{2}$
\\ \hline
$ D^{4}_{4} $ & $ a_{11}e_{1}\otimes e_{1}  + $ & $e_{1}\ast e_{1}= e_{1}$; $e_{1}\ast e_{2}= e_{2};$ $ e_{2}\ast e_{1}= e_{2};  e_{1}\ast e^{\ast}_{1}= e^{\ast}_{1} $ \cr
  & $ a_{22}e_{2}\otimes e_{2}   $   &  $ e^{\ast}_{2}\ast e^{\ast}_{2}= -a_{22}e^{\ast}_{1}; e^{\ast}_{1}\ast e^{\ast}_{2}= -a_{11}e^{\ast}_{2};$ $e^{\ast}_{1}\ast e^{\ast}_{1}= -a_{11} e^{\ast}_{1} $ \cr
  & &   $ e_{2}\ast e^{\ast}_{1}= -a_{11}e_{2}; e^{\ast}_{1}\ast e_{1}= -a_{11}e_{1}; e^{\ast}_{1}\ast e_{2}= -a_{11}e_{2} $\cr
  && $ e^{\ast}_{2}\ast e^{\ast}_{2}= a_{11}e^{\ast}_{1} -a_{11}e^{\ast}_{2}; e^{\ast}_{2}\ast e^{\ast}_{2}= a_{11}e_{1} -a_{22}e_{2} + e^{\ast}_{1}$ \cr
  && $ e^{\ast}_{2}\ast e_{1}= -a_{11}e^{\ast}_{1} + e^{\ast}_{2} - e^{\ast}_{1}; e^{\ast}_{2}\ast e_{2}= a_{11}e_{1} + e^{\ast}_{1}$
 \\ \hline
 $ D^{5}_{1} $ & $ a_{11}e_{1}\otimes e_{1} +  $ & $ e_{i}\ast e_{j}= e^{\ast}_{i}\ast e_{j}= 0$ \cr
$ \beta=0, 1 $  & $ a_{12}e_{1}\otimes e_{2} + $  &  $ e_{i}\ast e^{\ast}_{j}= e^{\ast}_{i}\ast e^{\ast}_{j}= 0 $ \cr
  & $ a_{12}e_{21}\otimes e_{1} + $   & $i, j= 1, 2$; $ \beta= 0 $    \cr
  & $ a_{22}e_{2}\otimes e_{2};   $  &  \cr
\cline{2-3}
 & $ a_{11}e_{1}\otimes e_{1}   $ & $ e^{\ast}_{1}\ast e_{2}=e^{\ast}_{2};$ $ e_{2}\ast e^{\ast}_{1}= -e^{\ast}_{2};$ $ \beta= 1 $

\\ \hline
 $ D^{6}_{1} $ & $ a_{11}e_{1}\otimes e_{1} +  $ & $ e_{2}\ast e_{2}= e_{2}$; $ e^{\ast}_{2}\ast e^{\ast}_{2}= -a_{22}e^{\ast}_{2}$ \cr
  & $ a_{22}e_{2}\otimes e_{2} $  &  $ e_{2}\ast e^{\ast}_{2}= -a_{22}e_{2}; e^{\ast}_{2}\ast e_{2}= e^{\ast}_{2} $
  \\ \hline
 $ D^{6}_{2} $ & $ a_{22}e_{2}\otimes e_{2}  $ & $ e_{2}\ast e_{2}= e_{2}; e^{\ast}_{2}\ast e_{2}= e^{\ast}_{2}$; $ e^{\ast}_{2}\ast e^{\ast}_{2}= - a_{22}e^{\ast}_{2}; e^{\ast}_{1}\ast e_{2}= e^{\ast}_{1}$  \cr
  &    &  $ e^{\ast}_{1}\ast e^{\ast}_{1}= a_{22} e^{\ast}_{1};$ $ e^{\ast}_{1}\ast e_{1}= a_{22} e_{1};$ $ e^{\ast}_{2}\ast e_{1}= -a_{22} e_{1}$\cr
  &    &  $ e^{\ast}_{2}\ast e^{\ast}_{2}= -a_{22} e_{2};$ $ e_{1}\ast e^{\ast}_{1}= -a_{22} e_{2} -e^{\ast}_{2}$ \cr
  \cline{2-3}
  & $ a_{11}e_{1}\otimes e_{1} + $ & $ e_{2}\ast e_{2}= e_{2}; e^{\ast}_{2}\ast e_{2}= e^{\ast}_{2}$; $ e^{\ast}_{1}\ast e_{2}= e^{\ast}_{1}; e_{1}\ast e^{\ast}_{1}= -e^{\ast}_{2}$   \cr
  & $ a_{12}e_{1}\otimes e_{2} + $ & $  e_{2}\ast e^{\ast}_{1}= -a_{12}e_{2}; e^{\ast}_{1}\ast e_{1}= a_{11}e_{1};  e^{\ast}_{2}\ast e_{1}= a_{12}e_{1} $ \cr
  & $ a_{12}e_{2}\otimes e_{1}  $ & $e^{\ast}_{1}\ast e^{\ast}_{2}= -a_{12}e^{\ast}_{2} + a_{12}e^{\ast}_{1}; e^{\ast}_{1}\ast e^{\ast}_{1}= a_{11}e^{\ast}_{1}$
  \\ \hline
  $ D^{6}_{3} $ & $ a_{22}e_{2}\otimes e_{2}  $ & $ e_{2}\ast e_{2}= e_{2}; e_{2}\ast e^{\ast}_{1}= e^{\ast}_{1}$; $ e_{2}\ast e^{\ast}_{2}= e^{\ast}_{2}$  \cr
  &    &  $ e^{\ast}_{2}\ast e^{\ast}_{1}= -a_{22} e^{\ast}_{1};$ $ e_{1}\ast e^{\ast}_{1}= - e^{\ast}_{2} -a_{22} e_{2}$ \cr
\cline{2-3}
  & $ a_{11}e_{1}\otimes e_{1} + $ &  $ e_{2}\ast e_{2}= e_{2}; e^{\ast}_{2}\ast e^{\ast}_{2}= a_{12} e^{\ast}_{2}$; $ e_{2}\ast e^{\ast}_{1}= -a_{12}e_{2}$ \cr
  & $ a_{12}e_{1}\otimes e_{2} + $ & $e_{2}\ast e^{\ast}_{1}= e^{\ast}_{1} + a_{11}e_{1};  e_{2}\ast e^{\ast}_{2}= e^{\ast}_{2} + a_{12}e_{2}$ \cr
  & $ a_{12}e_{2}\otimes e_{1}  $ & $e_{1}\ast e^{\ast}_{1}= -e^{\ast}_{2} - a_{12}e_{1};  e^{\ast}_{1}\ast e^{\ast}_{1}= a_{11}e^{\ast}_{2} - a_{12}e^{\ast}_{1}$
  \\ \hline
 $ D^{6}_{4} $ & $ a_{22}e_{2}\otimes e_{2} $ & $ e_{2}\ast e_{2}= e_{2}; e^{\ast}_{2}\ast e_{1}=  e^{\ast}_{1}$; $ e^{\ast}_{2}\ast e^{\ast}_{2}= -a_{22}e^{\ast}_{2}$ \cr
  &    &  $ e_{2}\ast e^{\ast}_{2}= -a_{22}e_{2};$ $ e_{1}\ast e^{\ast}_{2}= - e^{\ast}_{1}; e^{\ast}_{2}\ast e_{2}=  e^{\ast}_{2} $
\\ \hline
 $ D^{7}_{1} $ & $ a_{11}e_{1}\otimes e_{1} + $ & $ e_{1}\ast e_{1}= e_{1}$; $ e_{2}\ast e_{2}=  e_{2};$ $ e^{\ast}_{2}\ast e^{\ast}_{1}= -a_{12}e^{\ast}_{2}; e^{\ast}_{2}\ast e^{\ast}_{2}= -a_{22}e^{\ast}_{2}$ \cr
  &  $ a_{12}e_{1}\otimes e_{2} + $  &  $ e_{1}\ast e^{\ast}_{2}= -a_{11}e_{1}; e_{1}\ast e^{\ast}_{2}= -a_{12}e_{1};  e_{2}\ast e^{\ast}_{1}= -a_{12}e^{\ast}_{2}$  \cr
  & $ a_{21}e_{2}\otimes e_{1}  $ & $ e_{2}\ast e^{\ast}_{2}= -a_{22}e_{2};$ $ e^{\ast}_{1}\ast e_{2}= -a_{12}e_{2}; e^{\ast}_{2}\ast e_{1}= -a_{12}e_{1}$ \cr
  & $ a_{22}e_{2}\otimes e_{2}$ &  $ e^{\ast}_{1}\ast e^{\ast}_{1}= -a_{11}e^{\ast}_{1}; e^{\ast}_{1}\ast e^{\ast}_{2}= -a_{12}e^{\ast}_{1}  $                                                 \cr
  & $ a_{12}= a_{21} $ & $ e^{\ast}_{1}\ast e_{1}= a_{12}e_{2} + e^{\ast}_{1};$  $ e^{\ast}_{2}\ast e_{2}= a_{12}e_{1} + e^{\ast}_{2} $
\\ \hline
     $ D^{7}_{2} $ & $ a_{11}e_{1}\otimes e_{1} + $ & $ e_{1}\ast e_{1}= e_{1}$; $ e_{2}\ast e_{2}=  e_{2};$ $ e_{1}\ast e^{\ast}_{1}= -a_{11}e_{1}; e_{1}\ast e^{\ast}_{2}= -a_{12}e_{1} $ \cr
  &  $ a_{12}e_{1}\otimes e_{2} + $  & $ e^{\ast}_{1}\ast e_{2}= -a_{12}e_{2}; e^{\ast}_{2}\ast e_{2}= -a_{22}e_{2};$  $ e^{\ast}_{1}\ast e^{\ast}_{1}= (-a_{11}-a_{12})e^{\ast}_{1} $   \cr
  & $ a_{21}e_{2}\otimes e_{1}  $ &  $ e^{\ast}_{2}\ast e^{\ast}_{1}= -a_{22}e^{\ast}_{1}; e^{\ast}_{2}\ast e^{\ast}_{2}= -a_{22}e^{\ast}_{2};$ $ e^{\ast}_{1}\ast e^{\ast}_{2}= -a_{12}e^{\ast}_{2} -a_{12}e^{\ast}_{1} $  \cr
  & $ a_{22}e_{2}\otimes e_{2}  $ & $e_{2}\ast e^{\ast}_{1}= -a_{12}e_{2} + e^{\ast}_{1};$ $ e_{2}\ast e^{\ast}_{2}= a_{12}e_{1} + e^{\ast}_{2} $                                                 \cr
  & $ a_{12}= a_{21} $ & $ e^{\ast}_{1}\ast e_{1}= -a_{12}e_{1} +(a_{12}-a_{22})e_{2}  $
 $ + e^{\ast}_{1} - e^{\ast}_{2} $
 \\ \hline
      $ D^{7}_{3} $ & $ a_{11}e_{1}\otimes e_{1} + $  & $ e_{1}\ast e_{1}= e_{1}$; $ e_{2}\ast e_{2}= e_{2};$ $ e^{\ast}_{1}\ast e^{\ast}_{1}= -a_{11}e^{\ast}_{1}; e^{\ast}_{2}\ast e^{\ast}_{2}= -a_{22}e^{\ast}_{2} $ \cr
  & $ a_{22}e_{2}\otimes e_{2} $   &  $ e_{1}\ast e^{\ast}_{1}= -a_{11}e_{1}; e_{2}\ast e^{\ast}_{2}=  e^{\ast}_{2};$  $ e^{\ast}_{1}\ast e_{1}= e^{\ast}_{1}; e^{\ast}_{2}\ast e_{2}= -a_{22}e_{2} $
\\ \hline
\end{tabular}

\begin{center}
\begin{tabular}{c c c}
\hline
Dendriform   & Symmetric  &  Connes cocycles \cr
 algebra   & solutions     &  structures over $ \mathcal{A}\oplus  \mathcal{A}^{\ast} $\cr
structures &&
    \\ \hline
 $ D^{7}_{4} $ & $ a_{11}e_{1}\otimes e_{1} + $  & $ e_{1}\ast e_{1}= e_{1}$; $ e_{2}\ast e_{2}= e_{2};$  $ e^{\ast}_{2}\ast e^{\ast}_{1}= -a_{11}e^{\ast}_{2}; e^{\ast}_{1}\ast e^{\ast}_{2}= -a_{12}e^{\ast}_{2} $ \cr
  & $ a_{12}e_{1}\otimes e_{2} $   &  $ e_{1}\ast e^{\ast}_{1}= -a_{11}e_{1}; e_{1}\ast e^{\ast}_{2}=  -a_{12}e_{1}; e_{2}\ast e^{\ast}_{1}= -a_{12}e^{\ast}_{1}$  \cr
   &  $ a_{21}e_{2}\otimes e_{1} $  &  $e^{\ast}_{1}\ast e_{2}= -a_{12}e_{2};$ $ e^{\ast}_{2}\ast e_{1}= e^{\ast}_{2}; e^{\ast}_{1}\ast e^{\ast}_{1}= -a_{11}e^{\ast}_{1} $ \cr
   & $ a_{21}= a_{12}  $ & $e^{\ast}_{2}\ast e^{\ast}_{2}= -a_{12}e^{\ast}_{2};$ $ e^{\ast}_{1}\ast e_{1}= a_{12}e_{2} + e^{\ast}_{1} $ \cr
&&$ e_{2}\ast e^{\ast}_{2}= (-a_{11}+ a_{12})e_{1} -a_{12}e_{2} -e^{\ast}_{1} +e^{\ast}_{2}$
\\ \hline
  $ D^{7}_{5, \lambda} $  & $  a_{11}e_{1}\otimes e_{1}  $ & $ e_{1}\ast e_{1}= e_{1}$; $ e_{2}\ast e_{2}= e_{2};$ $ e_{1}\ast e^{\ast}_{1}= e^{\ast}_{1}; e^{\ast}_{1}\ast e_{1}=  -a_{11}e_{1} $ \cr
 && $ e_{1}\ast e^{\ast}_{2}= -\lambda e^{\ast}_{2}; e^{\ast}_{1}\ast e^{\ast}_{1}= -a_{11}e^{\ast}_{1};$ $ e^{\ast}_{1}\ast e^{\ast}_{2}= a_{11}\lambda e^{\ast}_{2} $ \cr
 &&$ e^{\ast}_{2}\ast e_{2}= a_{11}\lambda e_{1} +\lambda e^{\ast}_{1} + e^{\ast}_{2} $ \cr
\cline{2-3}
 & $ a_{22}e_{2}\otimes e_{2} $ & $ e_{1}\ast e_{1}= e_{1}$; $ e_{2}\ast e_{2}= e_{2};$ $ e_{1}\ast e^{\ast}_{1}= e^{\ast}_{1}; e_{1}\ast e^{\ast}_{2}= -a_{22}\lambda e_{2} $ \cr
   &    &  $ e_{2}\ast e^{\ast}_{2}= -a_{22} e_{2}; e_{1}\ast e^{\ast}_{2}= -\lambda e^{\ast}_{2};$   $ e^{\ast}_{2}\ast e_{2}= \lambda e^{\ast}_{1} + e^{\ast}_{2} $ \cr
  && $ e^{\ast}_{2}\ast e^{\ast}_{2}= -\lambda a_{22}e^{\ast}_{1} -a_{22}e^{\ast}_{2}  $ \cr
  \cline{2-3}
   & $ a_{11}e_{1}\otimes e_{1} $ & $ e_{1}\ast e_{1}= e_{1}$; $ e_{2}\ast e_{2}= e_{2};  e_{1}\ast e^{\ast}_{1}= e^{\ast}_{1};  e_{2}\ast e^{\ast}_{1}= -a_{12}e_{2}$ \cr
  & $ a_{12}e_{1}\otimes e_{2} $   &  $ e_{2}\ast e^{\ast}_{2}= -a_{22}e_{2};$ $ e^{\ast}_{1}\ast e_{1}= -a_{11} e_{1}; e^{\ast}_{1}\ast e_{2}= -a_{12}e_{2} $ \cr
   & $ a_{21}e_{2}\otimes e_{1} $   & $ e^{\ast}_{2}\ast e_{1}= -a_{12}e_{1}; e^{\ast}_{1}\ast e^{\ast}_{2}= -a_{12}e^{\ast}_{2};  e^{\ast}_{1}\ast e^{\ast}_{1}= -a_{11}e^{\ast}_{1}$  \cr
    & $ a_{22}e_{2}\otimes e_{2} $ & $e^{\ast}_{2}\ast e^{\ast}_{1}= -a_{12}e^{\ast}_{2};$  $ e^{\ast}_{2}\ast e_{2}= -a_{12}e_{2} -e^{\ast}_{1} + e^{\ast}_{2}  $  \cr
  &$ a_{21}= a_{12} $ & $ e^{\ast}_{2}\ast e^{\ast}_{2}= -a_{22}e^{\ast}_{1} + (-a_{22}-a_{12})e^{\ast}_{2}   $
  \\ \hline
 $ D^{7}_{6} $ & $ a_{11}e_{1}\otimes e_{1} $ & $ e_{1}\ast e_{1}= e_{1}$; $ e_{2}\ast e_{2}= e_{2};$ $ e_{1}\ast e^{\ast}_{1}= -a_{12}e_{1}+ e^{\ast}_{1} -e^{\ast}_{2} $ \cr
  & $ a_{12}e_{1}\otimes e_{2} $ & $ e_{1}\ast e^{\ast}_{2}= -a_{12}e_{1}; e_{2}\ast e^{\ast}_{1}= -a_{12}e_{2}; e_{2}\ast e^{\ast}_{2}= -a_{22}e_{2}$   \cr
   & $ a_{21}e_{2}\otimes e_{1} $   & $ e^{\ast}_{1}\ast e^{\ast}_{2}= -a_{12}e^{\ast}_{1} $; $ e^{\ast}_{2}\ast e_{1}= -a_{12}e_{1};$ $ e^{\ast}_{1}\ast e_{1}= -a_{11}e_{1}$  \cr
    & $ a_{22}e_{2}\otimes e_{2} $   &  $ e^{\ast}_{2}\ast e^{\ast}_{2}= (-a_{22} -a_{12})e^{\ast}_{2}; e^{\ast}_{1}\ast e_{2}= a_{11}e_{1}+ a_{12}e_{2} + e^{\ast}_{1} $  \cr
  &$ a_{21}= a_{12}=a_{22} $ & $ e^{\ast}_{2}\ast e^{\ast}_{1}= -a_{12}e^{\ast}_{1}  + a_{12} e^{\ast}_{2};$ $ e^{\ast}_{2}\ast e_{2}= a_{12}e_{1}  -a_{12}e_{2} + e^{\ast}_{2} $ \cr
 && $ e^{\ast}_{1}\ast e^{\ast}_{1}= (-a_{11}-a_{12})e^{\ast}_{1}+ a_{11}e^{\ast}_{2}  $
  \\ \hline
 $ D^{7}_{7} $ & $ a_{11}e_{1}\otimes e_{1} $ & $ e_{1}\ast e_{1}= e_{1}$; $ e_{2}\ast e_{2}= e_{2};$   $ e_{2}\ast e^{\ast}_{1}= e^{\ast}_{2}; e_{2}\ast e^{\ast}_{2}= e^{\ast}_{2} $ \cr
   &    &  $ e^{\ast}_{1}\ast e_{1}= -a_{11} e_{1};$ $ e^{\ast}_{1}\ast e_{2}= a_{11}e_{1}+ e^{\ast}_{1} -e^{\ast}_{2}  $\cr
    &    &  $ e_{1}\ast e^{\ast}_{1}= e^{\ast}_{1} -e^{\ast}_{2};$ $ e^{\ast}_{1}\ast e^{\ast}_{1}= -a_{11}e^{\ast}_{1} +a_{11}e^{\ast}_{2}  $
 \\ \hline
  $ D^{7}_{8, \lambda} $ & $ a_{22}e_{2}\otimes e_{2} $ & $ e_{1}\ast e_{1}= e_{1}$; $ e_{2}\ast e_{2}= e_{2};$ $ e_{1}\ast e^{\ast}_{2}= \lambda e^{\ast}_{1}; e_{2}\ast e^{\ast}_{2}= -a_{22}e_{2} $ \cr
   &    &  $ e^{\ast}_{1}\ast e_{1}= e^{\ast}_{1};  e^{\ast}_{2}\ast e_{1}= -\lambda e^{\ast}_{1};$  $ e^{\ast}_{2}\ast e_{2}= e^{\ast}_{2};  e^{\ast}_{2}\ast e^{\ast}_{2}= -a_{22}e^{\ast}_{2}$
     \\ \hline
  $ D^{7}_{9, \lambda} $ & $ a_{22}e_{2}\otimes e_{2} $ & $ e_{1}\ast e_{1}= e_{1}$; $ e_{2}\ast e_{2}= e_{2};$ $ e_{1}\ast e^{\ast}_{2}= \lambda e^{\ast}_{1}; e^{\ast}_{2}\ast e_{2}= -a_{22}e_{2} $ \cr
   &    &  $e_{2}\ast e^{\ast}_{2}= e^{\ast}_{2}; e^{\ast}_{1}\ast e_{1}= e^{\ast}_{1};  e^{\ast}_{2}\ast e_{1}= -\lambda e^{\ast}_{1};$  $e^{\ast}_{2}\ast e^{\ast}_{2}= -a_{22}e^{\ast}_{2}$
     \\ \hline
 $ D^{7}_{10, \lambda} $ & $ a_{22}e_{2}\otimes e_{2} $ & $ e_{1}\ast e_{1}= e_{1}$; $ e_{2}\ast e_{2}= e_{2};$  $ e_{1}\ast e^{\ast}_{2}= -a_{22} e_{2}- \lambda e^{\ast}_{1} - e^{\ast}_{2} $
 \cr
  &    &  $ e_{1}\ast e^{\ast}_{1}= e^{\ast}_{1}; e^{\ast}_{2}\ast e^{\ast}_{2}= -a_{22}e^{\ast}_{2};$   $ e_{2}\ast e^{\ast}_{2}= -a_{22} e_{2}- e^{\ast}_{1} $ \cr
    &    &  $ e^{\ast}_{2}\ast e_{2}= e^{\ast}_{1} + e^{\ast}_{2};$ $ e^{\ast}_{2}\ast e_{1}= \lambda e^{\ast}_{1} + e^{\ast}_{2}  $
   \\ \hline
 $ D^{7}_{11, \lambda} $ & $ a_{22}e_{2}\otimes e_{2} $ & $ e_{1}\ast e_{1}= e_{1}$; $ e_{2}\ast e_{2}= e_{2};$ $ e_{1}\ast e^{\ast}_{2}= -a_{22} e_{2} -\lambda e^{\ast}_{1} - e^{\ast}_{2} $ \cr
  &    &  $ e_{1}\ast e^{\ast}_{1}= e^{\ast}_{1}; e^{\ast}_{2}\ast e^{\ast}_{2}= -a_{22}e^{\ast}_{2} ;$   $ e_{2}\ast e^{\ast}_{2}=  e^{\ast}_{2}- e^{\ast}_{1} $ \cr
    &    &  $ e^{\ast}_{2}\ast e_{2}= -a_{22}e_{2} + e^{\ast}_{1};$  $ e^{\ast}_{2}\ast e_{1}= a_{22}e_{2} + \lambda e^{\ast}_{1} + e^{\ast}_{2}  $
\\ \hline
\end{tabular}
\end{center}

 \section{Concluding remarks}

 In this work, we gave an overview of the main concepts, definitions and known fundamental results related to the notions of Frobenius algebras, bialgebras, and Connes cocycles.
We classified the solutions of the associative Yang-Baxter equation on complex associative algebras in dimensions $1$ and $ 2 $. In dimension 1, we obtained $ 2 $ classes against $ 7 $ in dimension 2. The skew-symmetric solutions enabled us to carry out double constructions of the Frobenius algebras of these associative algebras.
 Finally, we obtained all compatible dendriform algebras  in dimensions $1$ and $ 2 $, gave a classification of solutions of D-equations and the  skew-symmetric Connes cocycles associated with  symmetric solutions.

 \section*{Aknowledgement}
 MNH is grateful to Prof C. Bai from Chern Institute of Mathematics, China, for useful discussions, inputs and provided references.
 This work is partially supported by the Abdus Salam International Centre for Theoretical
 Physics (ICTP, Trieste, Italy) through the Office of External Activities (OEA) - Prj-15. The
 ICMPA is also in partnership with the Daniel Iagolnitzer Foundation (DIF), France.

\end{document}